\documentclass{amsart}
\usepackage{amsfonts,amssymb, amsmath}

\numberwithin{equation}{section}

\newcommand{\RR}{\mathbb{R}}
\newcommand{\ZZ}{\mathbb{Z}}
\newcommand{\NN}{\mathbb{N}}
\newcommand{\TT}{\mathbb{T}}

\newcommand{\p}{\partial_{i}}
\newcommand{\pp}{\partial_{i+1}}
\newcommand{\pn}{\partial_{i-1}}

\newcommand{\lims}{\limsup_{\delta \to 0, N \to \infty}}
\newcommand{\limslN}{\limsup_{l \to \infty, N \to \infty}}
\newcommand{\e}{\epsilon}
\newcommand{\fri}{{\tiny \frac{i}{N}}}
\newcommand{\La}{\Lambda}

\newcommand{\Ave}[3]{\mathrm{Av}_{j= #2 - #3 }^{#2 + #3} \tau^{j} #1} 
\newcommand{\Av}[3]{\mathrm{Av}_{j= #2 - #3 }^{#2 + #3} #1_{j}}
\newcommand{\Avph}{\mathrm{Av}_{j=i-l}^{i+l} \phi(\frac{j}{N})}

\newcommand{\lr}[1]{ \left \langle \langle #1 \right \rangle \rangle }
\newcommand{\lrp}[2]{\left \langle \langle #1, #2 \right \rangle \rangle}
\newcommand{\VV}{\mathbb{V}}

\newtheorem{theorem}{Theorem}[section]
\newtheorem{definition}{Definition}[section]
\newtheorem{lemma}{Lemma}[section]

\begin{document}

\title{Hydrodynamic scaling limit of continuum solid-on-solid model}

\author{BY Anamaria Savu}
\address{Department of Mathematics, Queen's University, 
          Kingston, ON, K7L 3N6, Canada}
\email{ana@mast.queensu.ca}

\thanks{Research supported by an Ontario Graduate Studies Fellowship and the 
        University of Toronto Fellowship}





\begin{abstract}
A fourth order nonlinear evolution equation is derived from 
a microscopic model for surface diffusion, namely, the continuum 
solid-on-solid model. We use the method developed by Varadhan for the 
computation of the hydrodynamic scaling limit of nongradient
models. What distinguishes our model from other models discussed 
so far is the presence of two conservation laws for the dynamics in 
a nonperiodic box and the complex dynamics that is not 
nearest-neighbor. Along the way, a few steps has to be adapted to our new 
context. As a byproduct of our main result we also derive the
hydrodynamic scaling limit of a perturbation of continuum
solid-on-solid model, a model that incorporates both surface diffusion and surface electromigration.
\end{abstract}

\maketitle

\tableofcontents

\newpage


\section{Introduction}

A process of great technological importance, molecular beam epitaxy (MBE), 
is used to manufacture computer chips and  semiconductor
devices, see Barabasi and Stanley \cite{Bar}. In general the chip is 
constructed by spraying a beam of atoms on a flat surface.
There are three phenomena that take place in the construction 
process namely deposition, diffusion  and desorption. The atoms arrive or 
deposit on the surface and do not stick on the first contact point, 
but diffuse or walk on the surface. 
When an atom reaches the edge of another wandering atom, the two atoms meet 
or glue together, forming islands. Sometimes an atom may jump out of the surface. 
Smaller islands may develop into larger islands affecting the roughness of 
the surface on the macroscopic scale. A rough surface does not have very good 
contact properties and engineers would like to understand the basics
mechanisms affecting the morphology in general. 
 
How deposition and desorption affect the morphology of the surface, it is 
quite well understood. It is of great importance to know how the profile of 
the surface evolves on the macroscopic scale if the atoms that make up the 
surface diffuse. We assume that no atoms arrive or leave the surface.

The present paper discusses a model for surface diffusion the continuum 
solid-on-solid model, known also as the fourth-order Ginzburg-Landau model.
The system has a complex interaction that is not a nearest-neighbor 
interaction and has two conserved quantities in a nonperiodic box, 
namely the total slope of the surface and the linear mean of the surface slope.
Under the assumption that the molecules of the surface follow the dynamics 
of the continuum solid-on-solid model, we will prove that the dynamics of
the  surface slope profile on the macroscale is a fourth-order nonlinear equation.
Since the model is nongradient, the derivation of the limit is 
not trivial and we shall use  the method developed in Varadhan \cite{Var1} 
and further extended in Quastel \cite{Qua1}, Varadhan and Yau \cite{Var2}. 
Even though we have an interaction
that is not of nearest-neighbor type, the microscopic current still 
decomposes as the Laplacian of the slope field and the fluctuations. 
Usually in nearest-neighbor models after subtracting the fluctuations
from the current we end up with the gradient of some local function.
We use the result that for continuum solid-on-solid model the 
space of exact functions has codimension one inside the space of 
closed functions, see Savu \cite{Sav}. 

We also include a discussion of a perturbation of continuum
solid-on-solid model, where the evolution of the surface is driven by both 
diffusion and electric field. The electric field will add one extra 
second-order term to the final nonlinear equation.  

The continuum solid-on-solid model is an approximation of the discrete
solid-on-solid model and hence is not considered a truly microscopic model.
Unfortunately at the time of writing of this paper, we don't have
the required techniques to solve the discrete solid-on-solid model.
However, partial rigorous results and numerical analysis
are available for the discrete case, see Krug, Doobs and Majaniemi \cite{Kru}.

Finally the paper is organized as follows: section 2 contains the description 
of continuum solid-on-solid model and of the model for surface electromigration, 
the statement of the main 
results and a note on similar models considered in the literature
so far; section 3 outlines the proof of the main result, the 
computation of the scaling limit of continuum solid-on-solid model;
section 4 shows how the final nonlinear equation is identified;
in section 5  we calculate the asymptotics of central limit variance, 
the main ingredient used in section 4 to prove that the microscopic
current can be replaced by a multiple of the field Laplacian;
in section 6 we calculate the scaling limit of the modified model to 
incorporate surface electromigration.

\section{The  models} \label{model}

The continuum solid-on-solid model, is a 1-dimensional lattice system 
with continuous order parameter, used to study the evolution of an 
interface. The model describes the movement of the interface at the 
mesoscopic level and hence is not considered a truly microscopic 
model that captures all the aspects of the particle motion, but 
it has the advantage of being more suitable for computations.\\

{\bf The height model.} Let $\mathbb{T}$ be the torus represented as the 
interval $[0,1]$ with $0$ and $1$ identified. For each positive integer $N$, 
there are $N$ scaled periodic lattice points located at sites $i/N$ in 
$\mathbb{T}$, $ i=1, \dots ,N $. We shall denote by $h_i (t)$ the height 
of the surface at the site $i/N$, at time $t$. Also because of the periodicity 
of the lattice points, $h_{N+1}(t)=h_1(t)$. The energy function $H_N(h)$ of 
a height configuration $h$ is chosen to be invariant under the global 
translation $h_i(t) \longmapsto h_i(t)+c$ and has the form
$$H_N(h)= \sum_{i=1}^N V(h_{i+1}-h_i).$$

In this paper we shall assume that the potential is quadratic $V(x)=x^2$. 
We shall require the evolution in  time of the surface to preserve the sum 
of the heights, to have as invariant measure the infinite 
mass measure $e^{-H_N(h)}dh$, and to be reversible. All these properties are 
satisfied by the solution of the stochastic differential system:\\
$$ dh_i(t)= -\frac{N^4}{2}(w_i-w_{i-1})dt + N^2(\sqrt{a_{i}}dB_{i}-\sqrt{a_{i-1}}dB_{i-1}), \quad  1\leq i \leq N .$$

\noindent
Above, $a_i(h)=a(h_i-h_{i-1}, h_{i+1}-h_i, h_{i+2}-h_{i+1})$ where $
a$ is a function with bounded continuous first derivatives 
that satisfies $ 0 < 1/a^* \leq a(x_{-1},x_0,x_1) \leq a^* < \infty$. 
The $n$ copies of the Brownian motion $B_i, \; i=1, \dots ,n$ are independent.
We define the instantaneous current $w_i(h)$ of particles over the 
bond $\{ i, i+1 \}$
\begin{eqnarray}
w_i(h) &  = & (\partial_{-1}a-2\partial_{0}a+\partial_{1}a)(h_i-h_{i-1}, h_{i+1}-h_i, h_{i+2}-h_{i+1}) +  \nonumber\\
       &     & + a_i(h) (V ' (h_i-h_{i-1})-2V ' (h_{i+1}-h_i)+ V ' (h_{i+2}-h_{i+1}) ).
\end{eqnarray}
Because the height process preserves the sum of the heights of the surface, 
this dynamics models surface diffusion. \\

{\bf The slope process.} A crucial property of the dynamics of the 
heights is the gauge property, namely the dynamics invariance under 
the action of the group $G$ of translation in the $(1, \dots , 1)$ direction, 
$$G = \{ T: \RR^N \to \RR^N \; | \; T(x_1, \dots , x_N) = (x_1+c, \dots, x_N+c),\; c\in \RR \}. $$ 
Hence there exists an induced dynamics on the quotient space 
$\mathbb{R}^N/G$ of equivalence classes. A representative of an equivalence 
class is the slope configuration $x_i(t) = h_{i+1}(t)-h_i(t)$, $1 \leq i \leq N-1$ 
and $x_N(t)=h_1(t)-h_N(t)$. Note that $\sum_{i=1}^N x_i(t)=0$. \\
As a function of the slope configuration of the surface, the energy becomes
$H_N(x)= \sum_{i=1}^N V(x_i)$.  

In the sequel, we shall study the slope process rather than the height process. 
The slope process is reversible and has as equilibrium distribution, the product 
probability measure $d\nu_N^{\mathrm{eq}} =e^{-H_N(x)}/Z^N dx$. Below we write 
down the stochastic differential system, the generator, and the Dirichlet form 
of the slope process:
\begin{equation} \label{x1}
 dx_i(t)= \frac{N^4}{2}(w_{i+1}-2w_i+w_{i-1})dt + N^2(\sqrt{a_{i-1}}dB_{i-1}-2\sqrt{a_{i}}dB_{i}+\sqrt{a_{i+1}}dB_{i+1}) , 
\end{equation}
\begin{equation} \label{gen}
 N^4L_N(f)=\frac{N^4}{2}\sum_{i=1}^N a_i(x)(\pp -2\p +\pn)^2 f +w_i(x) (\pp-2\p+\pn)f,  
\end{equation}
\begin{equation}
N^4D_N(f)= \frac{N^4}{2}\int_{\mathbb{R}^N} \sum_{i=1}^N a_i(x) (\pp f-2\p f -\pn f)^2(x) \frac{e^{-H_N(x)}}{Z^N} dx.
\end{equation}
The factor $N^4$ in the generator with the lattice spacing of $1/N$ 
represents the scaling of space and time. This scaling is needed to 
observe a nontrivial motion in the limit.

The diffusion (\ref{x1}) is driven in the direction of the linear 
vector fields $X_i=\pp -2\p +\pn$, $1 \leq i \leq N$, therefore is not 
ergodic in the whole space $\mathbb{R}^N$. It becomes ergodic when 
restricted to  the hyperplane $x_1+\dots +x_N = N\bar{x}$ of average 
slope $\bar{x}$.
The unique equilibrium probability measure of the dynamics restricted to 
the hyperplane is the conditional probability $e^{-H_N (x)}/Z^N dx$ 
given $x_1+\dots +x_N = N\bar{x}$. 

Instantaneously, the slope profile of the surface decreases at some 
site $i$ twice as much as increases at the adjacent sites 
$i-1$, $i+1$. We see that any update of the slope 
configuration affects the slopes at three sites. This type of interaction,
known as the three site interaction, is quite complex and has been very 
rarely studied, so far.
 
We should note an integration by parts property of the current $w_i$, for 
all bounded and smooth functions $f$:
\begin{equation}
E^{\mathrm{eq}}[w_i \cdot f] = E^{\mathrm{eq}}[ a_i \cdot (\pp -2\p +\pn)f].
\end{equation}
The expected value above is with respect to the equilibrium measure 
$d\nu_N^{\mathrm{eq}}$.\\

{\bf The dynamics in a nonperiodic box.} Although we are concerned with 
the study of the slope process defined on a periodic space, we make use of 
a similar slope process evolving in a box with nonperiodic boundary. Below 
we describe this new dynamics and its main properties.
Suppose $\La$ is a box with finitely many sites of the lattice $\mathbb{Z}$. 
The infinitesimal generators, 
\begin{equation} \label{genla}
 L_{\La}(f)=\frac{1}{2}\sum_{i \in \La, i+1, i-1 \in \La} a_i(x)
(\pp -2\p +\pn)^2 f +w_i(x) (\pp-2\p+\pn)f, 
\end{equation} 
respectively,
\begin{equation} \label{geninf}
 L_{\infty}(f)=\frac{1}{2}\sum_{i \in \mathbb{Z}} a_i(x)(\pp -2\p +\pn)^2 
f +w_i(x) (\pp-2\p+\pn)f, 
\end{equation}
produce two diffusion processes.

In the case of the first dynamics, generated by $L_{\La}$, 
there is no transport of particles over the boundary of the box, 
so the box does not have periodic boundary. If the box 
$\La = [-l,l] \cap \ZZ$ or $\La=[i-l,i+l] \cap \ZZ$ we use 
the shorter notation $L_l$, respectively, $L_{i,l}$ for the 
generator $L_{\La}$. The second dynamics generated by $L_{\infty}$ 
is a dynamics on the infinite lattice $\mathbb{Z}$. 
The dynamics $L_{i,l}$ preserves the average slope, $y_{i,l}^1=\frac{x_{i-l}+\dots +x_{i+l}}{2l+1}$, 
and the linear mean of the slopes, $y_{i,l}^2=\frac{(-l)x_{i-l}+\dots +lx_{i+l}}{l(l+1)}$, 
inside the box $[i-l,i+l] \cap \ZZ$. 
The second conserved quantity, $y_{i,l}^2$, should
be understood as a boundary condition that is preserved in time.
The linear mean slope, $y_{i,l}^2$, 
is conserved in time because we have a model for surface diffusion and 
the total height of the surface does not changes as time passes by.
We will use two different notation $\bar{x}_{i,l}$ and $y^1_{i,l}$
for the mean slope of the field in the box centered at $i$, of size $l$.
Also $y_{i,l}$ stands for the vector of the conserved quantities, $(y_{i,l}^1, y_{i,l}^2)$.
As a convention we will drop the subscript $i$, meaning the center of the box,
if the box is centered at the origin.\\

{\bf Equilibrium measures.} We proceed to describe next the equilibrium
measures of the dynamics that we have introduced.
For the dynamics $L_{i,l}$, the grand canonical measure is the product probability measure 
$\nu_{\alpha,i,l}^{gc} = \bigotimes_{j=i-l}^{i+l} \frac{e^{\alpha x_j-V(x_j)}}{Z(\alpha )}dx_j,$ 
whereas the canonical measure is the conditional probability measure 
$\nu_{y,i,l}^c= \nu_{\alpha,i,l}^{gc} (\quad | y_{i,l})$
 given  the level set 
$$ \{ x \in \RR^{2l+1} \; | \; \frac{x_{i-l}+\dots +x_{i+l}}{2l+1}
= y^1_{i,l}, \frac{(-l)x_{i-l}+\dots + lx_{i+l}}{l(l+1)}= y^2_{i,l} \}.$$ 
The canonical measure
is the unique stationary probability measure for the restricted dynamics on this set.
The Dirichlet forms of the operators $L_l$, $L_{i,l}$, with respect to the 
grand canonical measure with $\alpha=0$, are denoted by $D_l(f)$, 
respectively, $D_{i,l}(f)$. The Dirichlet form of the operator $L_l$ with
respect to the canonical measure is denoted by $D_{\nu^c_{y,l}}(f)$.
For $L_{\infty}$, the product measures $\nu_{\alpha}^{gc} = 
\bigotimes_{i \in \mathbb{Z}} \frac{e^{\alpha x_i-V(x_i)}}{Z(\alpha )}dx_i$ 
are equilibrium measures. We shall use the notation $\nu^{\mathrm{eq}}_N$ 
for the product probability measure $\bigotimes_{i=1}^{N} \frac{e^{-V(x_i)}}{Z}dx_i$.\\

{\bf  A model for surface electromigration.}
We consider also a perturbation of the continuum solid-on-solid model. 
The new system describes the evolution of a one-dimensional surface driven by
both surface diffusion and surface electromigration. The surface
electromigration refers to the motion of atoms on a solid surface 
that is caused by an electric current in the material. The electric
field interacts with the atoms of the surface as the wind blows the
sand particles, and a ripple pattern is observed in the long run. 
Electromigration along interfaces is believed to play a crucial role
in the failure of metallic circuits, see Schimschak and Krug \cite{Sch}, 
for further details.

We assume the electric field is a continuous function $E(t, \theta)$
defined on $[0,T] \times \TT$. As before $x_i$ represents the slope of the 
surface at the site $i/N$. 
The generator of the system on a periodic lattice, 
that incorporates the action of the electric field, is
\begin{equation} \label{elect1}
N^4L_{N,E}(f) = N^4L_N(f) + \frac{N^2}{2} \sum_{i=1}^N E\bigg(t, \frac{i}{N}\bigg)
             a(x_{i-1},x_i,x_{i+1})(\partial_{i-1}-2\partial_i+\partial_{i+1}) f. 
\end{equation}

{\bf Hydrodynamic scaling limit of the models.}
We shall call  $P_{N,T}^{\mathrm{neq}}$ and $P_{N,T}^{\mathrm{eq}}$ 
the law up to time $T$ of the slope process (\ref{x1}) started in some 
nonequilibrium distribution $\nu_N^{\mathrm{neq}}$, equilibrium 
distribution $\nu_N^{\mathrm{eq}}$, respectively. The law up to time 
$T$ of the slope process (\ref{elect1}), driven by the electric field, 
started in the nonequilibrium 
measure $\nu_N^{\mathrm{neq}}$ shall be called $P_{N,E,T}^{\mathrm{neq}}$.
 
Under $P_{N,T}^{\mathrm{neq}}$ and $P_{N,T}^{\mathrm{eq}}$ the random variable
\begin{equation} \label{emp}
\pi_N(t) = \frac{1}{N}\bigg( x_1(t) \delta_{\frac{1}{N}}+\dots + x_N(t) 
\delta_{\frac{N}{N}} \bigg)
\end{equation} 
has distributions $Q_{N,T}^{\mathrm{neq}}$ and $Q_{N,T}^{\mathrm{eq}}$, respectively.
We also refer to the random variable (\ref{emp}) as the empirical distribution. 
Every realization of this random variable is a measure-valued continuous path 
and $Q_{N,T}^{\mathrm{neq}}$ is a
distribution on the space 
$ \mathcal{X} = \bigcup_{(l\geq 0)} C([0, T], \mathcal{M}_l )$. 
The space $\mathcal{X}$ is endowed with the inductive limit topology, 
the strongest topology that makes all the inclusions 
of $C([0, T], \mathcal{M}_l ) $ continuous. The space of signed measures 
$\mathcal{M}_l$, with total variation not exceeding $l$, is a metrizable
space with the weak topology.

We say that a model has hydrodynamic scaling limit if under certain 
assumptions the sequence of laws of empirical distributions has a     
limit that is supported on the solution of an initial value problem. 
\begin{definition}
A sequence of initial distributions $\nu_N^{\mathrm{neq}}$ on 
$\mathbb{R}^N$ is said to correspond to the macroscopic slope 
profile $m_0 \in L^1(\mathbb{T})$ 
if the random variable $\pi_N$ converges weakly in probability to 
$\delta_{m_0(\theta) d\theta} $, i.e., for any continuous function 
$\phi \in C(\mathbb{T})$  and any $\e >0$,
$$\limsup_{N \to \infty} \nu_N^{\mathrm{neq}} \bigg\{ 
\bigg| \frac{1}{N} \sum_{i=1}^{N} \phi \bigg(\frac{i}{N}\bigg) x_i  - 
\int_{\mathbb{T}} \phi (\theta ) m_0 (\theta ) d \theta \bigg| > \e \bigg\} =0 .
$$
\end{definition}

Our first result of the paper says that the continuum 
solid-on solid model has a scaling limit and provide the form of the limiting 
evolution equation. 
More precisely,

\begin{theorem} \label{mainres}
Assume that the potential $V(x)$ is equal to $x^2/2$. 
Let $m_0 \in L^1(\mathbb{T})$ be a macroscopic slope profile such that 
$\int_{\mathbb{T}} m_0(\theta ) d \theta =0 $. 
Assume that the sequence of initial distributions 
$\{\nu_N^{\mathrm{neq}}\}_N$ corresponds to the profile $m_0$ 
and the initial relative entropy $H(\nu_N^{\mathrm{neq}} | \nu_N^{\mathrm{eq}})$ 
is of order $\mathcal{O}(N)$. 
Then the sequence of probability measures $\{ Q_{N,T}^{\text{neq}}\}_{N \geq 0}$ 
is tight in $\mathcal{X}$.
 
Any possible limit, $Q_T$, of a convergent subsequence of $\{ Q_{N,T}^{\text{neq}}\}_{N \geq 0}$ 
is concentrated on the weak solutions 
$(m(t,\theta)d \theta )_{t \in [0,T]}$ of the Cauchy problem with periodic boundary conditions
\begin{equation} \label{1eq1}
\partial_t m  =  - \frac{1}{2}\partial_{\theta}^2 ( \hat{a} (m) 
\partial_{\theta}^2 m ),  \quad \quad m(0, \theta)  =  m_0(\theta), \quad \theta \in \mathbb{T}.
\end{equation}
The transport coefficient $\hat{a}$ is a nonrandom continuous 
function on $\RR$, given by the following variational formula,
\begin{equation} \label{coeff}
\hat{a}(\alpha ) = \inf_g E_{\nu^{gc}_{\alpha}}\bigg[a(x_{-1},x_0,x_1)
\bigg(1+(\pp -2\p +\pn )\Big(\sum_{j \in \ZZ} \tau^{j}g\Big)\bigg)^2 \bigg].
\end{equation}
The infimum on the line above is taken over all local functions 
$g(x_{-s}, \dots ,x_s, \bar{y}^1_{0,l})$ of 
the slope configuration. The shift $\tau^j$ acts on configurations 
$(\tau^j x)_k = x_{k+j}$ and on local functions, 
$(\tau^j g)(x)= g(\tau^j x)$. The expectation $E_{\nu^{gc}_{\alpha}}$ is 
with respect to the grand canonical measure $\nu^{gc}_{\alpha}$. 
\end{theorem}

{\it Note.} It can be shown that the transport coefficient $\hat{a}$ is
a continuous, bounded above and below function (see Kipnis and Landim
\cite{Kip}). Also we expect that the methods of Landim, Olla and Varadhan 
\cite{Lan} can show that $\hat{a}$ is smooth but we will not pursue it here.
 
{\it Note.} By a weak solution of the Cauchy problem 
(\ref{1eq1}) we mean a path $m\in \mathcal{X}$  such that
for each time $T \geq 0$, the value of the path $m$ at the time $T$ is a 
Lebesgue absolutely continuous measures on $\mathbb{T}$, satisfying the
energy estimate
\begin{equation} \label{eneres}
\int_0^T \int_{\mathbb{T}} (\partial_{\theta}^2 m(t, \theta ))^2 d\theta dt < \infty .
\end{equation}
Moreover for each $0 \leq T < \infty $ and for each test function 
$\phi \in C^{1,2}([0,T] \times \mathbb{T})$, 
$$\int_{\mathbb{T}} m(T, \theta ) \phi (T, \theta) d \theta - 
\int_{\mathbb{T}} m(0, \theta ) \phi(0, \theta ) d \theta - 
\int_0^T \int_{\mathbb{T}} m(s, \theta ) \partial_s \phi (s, \theta ) d \theta ds +$$
$$+\frac{1}{2}\int_0^T \int_{\mathbb{T}} \hat{a} (m) \partial_{\theta}^2 m(s, \theta )
\partial_{\theta}^2 \phi(s, \theta ) d \theta ds =0. $$
If we assume that the initial condition $m_0$ has the property that 
$\int_{\TT } m_0(\theta ) d\theta =0$ then for each time $t \geq 0$ 
the solution satisfies $\int_{\TT } m(t, \theta ) d\theta =0$.

Uniqueness of weak solutions of the Cauchy problem (\ref{1eq1}), that satisfy
the energy estimate (\ref{eneres}) has not been proved yet. If the transport 
coefficient $\hat{a}$ does not depend on the field $m$, the uniqueness 
of the Cauchy problem is known and can be found in Eidelman book \cite{Ede}. 

As will be explained later the fluctuation dissipation equation for the continuum solid-on-solid
model follows from the direct sum decomposition of a Hilbert space to be defined next. 

We define the Hilbert space of closed functions, $\mathcal{C}_{X}$ to be 
the space of those $\xi \in L^2(d \nu_{\alpha}^{gc})$ that satisfy in the weak sense 
the equations $X_i(\tau^j \xi) =X_j(\tau^i \xi)$ for all integers $i$ and $j$.
It is not hard to see that a subspace  of $\mathcal{C}_{X}$ is the closed linear 
span in $L^2(d\nu_{\alpha}^{gc})$ of functions $\xi_g = X_0(\sum_{j \in \ZZ} \tau^j g)$, 
where $g$ is a bounded local function with bounded first derivatives. 
Even though the infinite sum $\sum_{j \in \ZZ} \tau^j g$ does not make sense, 
the function $\xi_g$ is well defined because the vector field kills 
all but finitely many terms of the infinite sum. We shall call this space the space of
exact functions and we shall use the notation $\mathcal{E}_X$. The space of exact
functions has codimension one inside the space of closed functions. 
In this paper we do not include the proof of this result,  since it is very technical 
and is not of probabilistic nature, however we include the statement, see Lemma \ref{clos1}.
 Lemma \ref{clos1} is discussed in  Savu \cite{Sav}.

\begin{lemma} \label{clos1}
Let ${\bf 1}$ denote the constant function $1$. The direct sum decomposition holds:
\begin{equation} \label{clos2}
\mathcal{C}_X = \RR {\bf 1} \bigoplus \mathcal{E}_X.  
\end{equation} 
\end{lemma}

The second result of the paper proves that the model for surface electromigration 
(\ref{elect1}) has a scaling limit as well, and calculates the limiting evolution equation.

\begin{theorem} \label{mainelect}
Suppose the hypothesis of Theorem \ref{mainres} are satisfied. Then 
the sequence of probability measures $\{ Q^{\mathrm{neq}}_{N,E,T}\}_{N \geq 0}$
is tight in $\mathcal{X}$, and any possible limit $Q_{E,T}$
is supported on the weak solutions of the Cauchy problem
\begin{equation}  \label{noneqel}
\partial_t m  =  - \frac{1}{2}\partial_{\theta}^2 ( \hat{a} (m) (\partial_{\theta}^2 m +E)), 
\quad
m(0, \theta)  =  m_0(\theta), \quad \theta \in \mathbb{T}. 
\end{equation}
The transport coefficient $\hat{a}$ is given by the same variational
formula as in the statement of Theorem \ref{mainres}.
\end{theorem}

{\bf Similar models.} 
The continuum solid-on-solid  model belongs to a large class 
of Ginzburg-Landau models. The slope model is of nongradient type 
and has an unusual dynamics because two neighboring exchanges occur 
always simultaneously. Nishikawa \cite{Nis}, and Bertini, Olla and Landim \cite{Ber} 
have investigated the hydrodynamic scaling limit  
in higher dimension of the gradient version (i.e., $a$ is a 
constant function) of the slope model. They have found that 
on the macroscale the interface follows a fourth-order nonlinear 
evolution equation.

Another similar model, the second-order Ginzburg-Landau model, where 
the sum of the heights is not conserved, 
was the subject of extensive discussions in the literature: 
the hydrodynamic scaling limit was derived by Fritz \cite{Fri} 
and  Guo, Papanicolau and Varadhan \cite{Guo} for the gradient version, and  
by Varadhan \cite{Var1} for the nongradient case, whereas 
the nonequilibrium fluctuations have been proved by Chang and 
Yau \cite{Cha}.
The second-order Ginzburg-Landau model and the continuum solid-on-solid 
model correspond to Glauber, respectively, Kawasaki dynamics in the context 
of interacting particle systems. 
As expected, a different dynamics at the mesoscopic level causes 
different dynamics at the macroscopic level, a second-order parabolic 
differential equation in the case of second-order Ginzburg-Landau model 
versus a fourth-order parabolic differential equation for the 
continuum solid-on-solid model.


\section{Hydrodynamic scaling limit of continuum solid-on-solid model}

In this section we give a sketch of the main result, Theorem \ref{mainres}.
We follow a standard scheme to derive the hydrodynamic scaling
limit of the continuum solid-on-solid model. The existence of
the limit follows from the tightness of the sequence of 
probability measures $\{Q_{N,T}^{\mathrm{neq}}\}_N$. Let $Q_T$ be
the limit of some weakly convergent subsequence of
$\{Q_{N,T}^{\mathrm{neq}}\}_N$. We proceed to characterize 
the limit $Q_T$, showing that it is supported on continuous paths
with certain regularity property, known as the energy estimate 
(\ref{eneres}). The most involved part of the argument is the 
identification of the possible weak limit $Q_T$ as some probability 
measure supported on the weak solution of the Cauchy problem (\ref{1eq1}).  
We shall make the assumption that the sequence of 
initial distributions $\{\nu_N^{\mathrm{neq}} \}_N $ 
corresponds to some macroscopic profile $m_0 \in L^1(\mathbb{T})$.

{\it Note on notations.} Throughout the paper we shall make use of the shorter 
notation $\limsup_{z_{1} \to i_1, \dots , z_n \to i_n} f (z_1, \dots , z_n)$ for the 
sequence of limits \\ $\limsup_{z_1 \to i_1} \dots \limsup_{z_n \to i_n} 
f(z_1, dots , z_n)$.\\

{\it Tightness.} 
As a consequence of Prohorov theorem and Arzela-Ascoli theorem, 
the tightness of the sequence $Q_{N,T}^{\mathrm{neq}}$ follows from
 the next two Lemmas. 

\begin{lemma} \label{lt1}
For any test function $\phi \in C^2(\mathbb{T})$, 
any finite time $T$, and any $\e >0$,
\begin{equation} \label{t2}
\lims P_{N,T}^{\mathrm{neq}} \bigg\{ \sup_{|s-t| \leq \delta, \; 0 \leq s,t \leq T} \bigg| 
\frac{1}{N} \sum_{i=1}^{N} \phi\bigg(\frac{i}{N}\bigg) x_i (t) - 
\frac{1}{N} \sum_{i=1}^{N} \phi\bigg(\frac{i}{N}\bigg) x_i(s) \bigg| > 
\e \bigg\} =0. 
\end{equation}
\begin{equation} \label{t3}
\limslN P_{N,T}^{\mathrm{neq}} \bigg\{ \sup_{0 \leq t \leq T} 
\frac{1}{N} \sum_{i=1}^{N} |x_i (t)| > l \bigg\} =0. 
\end{equation}
\end{lemma}
\noindent
{\bf Proof.} The first convergence (\ref{t2}) follows from 
Garsia-Rodemich-Rumsey inequality. To prove the second convergence 
(\ref{t3}) we can use estimates on the moment generating function 
of hitting time of the diffusion process (\ref{x1}). The reader 
may consult Kipnis and Landim \cite{Kip} or Guo, Papanicolau and 
Varadhan \cite{Guo} or Savu \cite{Sav} for a complete proof.

It is interesting to note that the stronger superexponential 
estimates can be established for the process in equilibrium,
\begin{equation} \label{t1}
\limsup_{\delta \to 0, N \to \infty} 
\frac{1}{N} \log P^{\mathrm{eq}}_{N,T} \bigg\{ \sup_{0 \leq t,s \leq T, t-s 
\leq \delta} \bigg| \int_s^t \frac{1}{N} \sum_{i=1}^N N^2 w_i(u) \phi\bigg(\fri\bigg) 
du \bigg| \geq \epsilon \bigg\} = - \infty. 
\end{equation}
\begin{equation} \label{t7}
\limslN \frac{1}{N} \log P_{N,T}^{\mathrm{eq}} \bigg\{ \sup_{0 \leq t \leq T} 
\frac{1}{N} \sum_{i=1}^{N} |x_i (t)| \geq l \bigg\} =- \infty .
\end{equation}
\qed

{\it Energy estimate.} Any limiting point $Q_T$ of the measure-valued 
sequence $\{ Q_{N,T}^{\mathrm{neq}} \}_N$ is supported on 
paths $\mu \in \mathcal{X}$  such that 
at each time $t$, $\mu (t)$ is a Lebesgue absolutely continuous 
measure on the torus $\TT$ with density $m(t, \theta)$.
Moreover for each finite time $T$, the density $m(t, \theta )$ 
satisfies the energy estimate
\begin{equation} \label{enerest}
\int_0^T  \int_{\mathbb{T}} (\partial_{\theta}^2 m(t, \theta ))^2 d\theta  dt < \infty .
\end{equation}
We note that the energy estimate ($\ref{enerest}$) is equivalent 
to the inequality
\begin{equation} \label{ener}
\sup_{\phi \in C^{1,2} ([0,T] \times \mathbb{T} ) } \int_0^T \int_{\mathbb{T}} 
(2 m (t, \theta ) \partial_{\theta}^2 \phi - C \phi^2 ) 
d\theta dt <  \infty,
\end{equation}
where $C$ is a constant not depending on $\phi$.

The entropy inequality (\ref{entropy}), Feynman-Kac formula (\ref{feynman}), and 
Lemma \ref{clt2} can be used to derive the estimate
$$ \limsup_{l,N \to  \infty } \sup_{\phi \in C^{1,2}([0,T] \times \mathbb{T})}E^{\mathrm{neq}} 
\bigg[ \frac{1}{N} \sum_{i=1}^N \int_0^T 
2 \phi\bigg(t, \fri \bigg) N^2\mathrm{Av}_{j=i-l+1}^{i+l-1}(\Delta x)_j  - $$
\begin{equation} \label{ener2}
-\VV_{i,l}(\Delta x, y) \phi^2\bigg(t, \fri \bigg) dt \bigg]< \infty .
\end{equation}
Here, the cylinder function $\Delta x$ is the discrete Laplacian of 
the slope field,  $x_{-1}-2x_0+x_1$, and the variance $\VV_{i,l}(\Delta x, y)$ 
is defined later in section 4, see definition \ref{}. 
As will be proved in Lemma \ref{boundonvar}, 
the variance, $\VV_{i,l}(\Delta x,y)$, has a uniform-in-$y$ lower bound, 
therefore the estimate (\ref{ener}) follows from (\ref{ener2}), 
after integrating by parts. Moreover, we can conclude that at any time
 $t$ the weak second derivative of the measure $\mu(t)$  is in 
$L^2(\mathbb{T},d\theta )$ and hence $\mu(t)$ is absolutely continuous 
with respect to Lebesgue measure on the torus $\mathbb{T}$. \\

{\it Identification of the equation.}
That any limiting point of the measure-valued sequence $\{Q_{N,T}^{\mathrm{neq}}\}_{N>0}$ 
is supported on the weak solutions of (\ref{1eq1})
follows if the event corresponding to the violation of the limiting equation 
has probability zero in the limit.

For each test function $\phi \in C^{1,2}([0,T] \times \mathbb{t})$ that 
are differentiable in time and twice differentiable in space, and each 
finite time $T$  we define the function
\begin{eqnarray} 
& & \quad \quad V(t) = \frac{1}{N} \sum_{i=1}^N \phi\bigg(t,\fri \bigg) x_i (t)  
    -\frac{1}{N} \sum_{i=1}^N \phi\bigg(0,\fri \bigg) x_i (0) - \\ \nonumber
& & \quad \quad \quad -\int_0^t \frac{1}{N} \sum_{i=1}^N \partial_s \phi (s, \fri)
    \bar{x}_{i,aN}(s) ds +\frac{1}{2}\int_{0}^{t} \frac{1}{N} 
    \sum_{i=1}^N \phi ''\bigg(s,\fri \bigg) \hat{a} (\bar{x}_{i,aN} (s)) \times \\ \nonumber 
& & \quad \quad \quad \quad \times b^{-2}(\bar{x}_{i-bN, cN} (s)-2\bar{x}_{i,cN}(s)
    +\bar{x}_{i+bN, cN}(s)) ds,
\end{eqnarray}
and the event 
\begin{equation} \label{event}
O^{\phi}_{a,b,c,\e}=\bigg\{ \sup_{0 \leq t \leq T} |V(t)|>\e \bigg\}.
\end{equation}
We will prove in the next sections that for each $\e >0$ we have
\begin{equation} \label{proba}
\limsup_{a,b,c \to 0, N \to \infty} P_{N,T}^{\mathrm{neq}} (O^{\phi}_{a,b,c,\e}) = 0.
\end{equation}
The proof of the result (\ref{proba}) is complicated and is divided into several steps. 
The function that defines this event  can be written as a 
sum of functions, see the beginning of section \ref{ident}. 
We will deal separately with each function in the sum and 
show that it converges to zero in probability. 


%
\section{Identification of the limiting equation} \label{ident}
In this section we establish that the event (\ref{event}) is negligible in
the limit, and hence any weak limit $Q_T$ is supported on the solutions 
of the Cauchy problem (\ref{1eq1}).
To save space, we are suppressing the time dependence of the test function $\phi$ 
that defines the event (\ref{event}). 

{\it Note on notations.} Assume $f$ is some local function. We denote by 
$\mathrm{Av}_{j=i-l}^{i+l} \tau^j f $ the average of shifts of $f$, namely 
$$\frac{\tau ^{i-l} f + \dots + \tau^{i+l}f}{2l+1}$$. 

We write $V(t)=$
\begin{eqnarray}
 &  & \\
 &  & =\frac{1}{N} \sum_{i=1}^N \phi\bigg(\fri \bigg) x_i(t) - \frac{1}{N} 
      \sum_{i=1}^N \phi\bigg(\fri \bigg) x_i (0) - \int_0^t \frac{1}{N} 
      \sum_{i=1}^N \phi\bigg(\fri \bigg) N^4 L_N( x_i)ds +  \nonumber \\
 &  & \\
 &  & \quad + \frac{1}{2}\int_0^t \frac{1}{N} \sum_{i=1}^N 
     \bigg\{ N^2 \bigg[ 
     \phi\bigg(\frac{i-1}{N}\bigg) -2 \phi \bigg(\fri\bigg) + 
     \phi \bigg(\frac{i+1}{N}\bigg) 
     \bigg]  -\phi '' \bigg(\fri\bigg) \bigg\} w_i  ds + \nonumber \\
 &  &\quad \quad + \frac{1}{2}\int_0^t \frac{1}{N} \sum_{i=1}^N \phi''
      \bigg(\fri \bigg)N^2 (w_i - \mathrm{Av}_{j=i-l_1}^{i+l_1} w_j) ds + \\
 &  &\quad \quad + \int_0^t \frac{1}{N} \sum_{i=1}^N \phi''\bigg(\fri \bigg) N^2 
      \bigg[ \mathrm{Av}_{j=i-l_1}^{i+l_1} 
      w_j - \hat{a}(\bar{x}_{i,l}) \Av{(\Delta x)}{i}{l_1} - \\
 &  & \quad \quad \quad \quad -\Ave{L_{\infty} f_r}{i}{l_1} \bigg] ds +  \nonumber \\
 &  & \quad \quad +  \frac{1}{2}\int_0^t \frac{1}{N}\sum_{i=1}^N \phi''\bigg(
     \fri \bigg) N^2 L_{\infty} (\Ave{f_r}{i}{l_1})  ds  \\
 &  & \quad \quad + \frac{1}{2}\int_0^t \frac{1}{N}\sum_{i=1}^N \phi''(\fri ) 
      N^2 \bigg[ \hat{a}(\bar{x}_{i,l}) \Av{(\Delta x)}{i}{l_1} - \\
 &  & \quad \quad \quad \quad -\hat{a} (\bar{x}_{i,aN}) b^{-2} (\bar{x}_{i-bN, cN}
     -2\bar{x}_{i,cN}+\bar{x}_{i+bN, cN}) \bigg] ds . \nonumber
\end{eqnarray}
We proceed to prove that each term in the sum above converges to $0$ 
in probability. \\

{\bf Martingale estimate (the term (4.1)).} 
We call $M_N(t)$ the term (4.1). From Ito formula we know that the 
process $\{M_N(t)\}_{t \geq 0}$ is a martingale and
$$E^{\mathrm{neq}}[M_N^2(T)] = E^{\mathrm{neq}}\bigg[ \int_0^T\frac{1}{N^2} 
\sum_{i=1}^N N^4 \bigg[ \phi\bigg(\frac{i-1}{N}\bigg) -2 \phi 
\bigg(\fri\bigg) + \phi \bigg(\frac{i+1}{N}\bigg)\bigg]^2 a_{i}ds \bigg].$$ 
The test function $\phi$ is chosen to have continuous second derivatives, 
therefore $E^{\mathrm{neq}}[M_N^2(T)]$ is of order $\mathcal{O}(\frac{1}{N})$.
We can use the Doob's inequality 
$$P^{\mathrm{neq}}_{N,T} \bigg\{ \sup_{0\leq t \leq T} | M_N(t) | \geq \epsilon\bigg\} 
\leq  \frac{1}{\epsilon^2} E[M_N^2(T)], \quad \epsilon > 0, $$
to conclude that the martingale is negligible in the limit, i.e.,
$$\lim_{N \to \infty} P^{\mathrm{neq}}_{N,T} \bigg\{ \sup_{0\leq t 
\leq T} | M_N(t) | \geq \epsilon \bigg\}  =0,  \quad \e > 0. $$

{\bf The term (4.2).}
A straightforward computation involving the Chebyshev inequality and 
the entropy inequality (\ref{entropy}) proves that the term (4.2) 
converges in probability to $0$. The test function needs to have
continuous fourth derivative. \\

{\bf A technical Lemma.}
We shall prove a Lemma that reduces the problem of establishing the 
negligibility of an event to finding that the largest eigenvalue 
of a Schr\"odinger operator is negative. For an operator 
$A : \mathcal{H} \to \mathcal{H}$ acting on a Hilbert space 
$\mathcal{H}$  we denote by $\mathrm{supspec}_{\mathcal{H}} A$ the 
largest value in the spectrum of $A$.    

\begin{lemma} \label{lem1}
Let $\{x(t)\}_{t \geq 0}$ be the slope  process with generator (\ref{gen}). 
Under the assumption 
\begin{equation}  
\label{sup1} 
\limsup_{N \to \infty} \mathrm{supspec}_{L^2(\nu_N^{\mathrm{eq}})} \bigg( \alpha g + \frac{N^4}{N} L_N \bigg) = 
\end{equation}
\begin{equation} \nonumber
=\limsup_{N \to \infty} \sup_{\rho, E^{\mathrm{eq}} [\rho^2]=1}  \bigg[ \alpha E^{\mathrm{eq}} [g \rho^2] - \frac{N^4}{N} D_N(\rho)\bigg]\leq 0, \quad  \alpha \neq 0
\end{equation}
it follows that the event $ \Big\{ \Big| \int_0^T g(x(s))ds \Big| \geq \epsilon \Big\}$ has negligible probability or:
\begin{equation} \label{neglig}
\lim_{N \to \infty} P^{\mathrm{neq}}_{N,T} \bigg\{ \bigg| \int_0^T g(x(t))dt \bigg| \geq \epsilon \bigg\} =0,  \; \; 
\; \epsilon >0.
\end{equation}
\end{lemma}

\noindent
{\bf Proof.}
We can use Chebyshev inequality to reduce the proof of (\ref{neglig}) to
\begin{equation} \label{l1}
\lim_{N \to \infty} E^{\mathrm{neq}} \bigg[ \bigg| \int_0^T g(x(t))dt \bigg| \bigg] =0  .
\end{equation}
Since we do not have much information about the initial 
nonequilibrium distribution we use the entropy inequality 
to replace the nonequilibrium distribution in (\ref{l1})
by the equilibrium distribution. 

Before we continue, we remind that given two probability measures 
$\nu$ and $\mu$ on the same probability space such that $\nu$ is absolutely continuous
with respect to $\mu$, we define the relative entropy of $\nu$
with respect to $\mu $ by 
$H(\nu|\mu)=E_{\mu}\Big[\frac{d\nu}{d\mu} \log \frac{d\nu}{d\mu} \Big]$ 
where $\frac{d\nu}{d\mu}$ is the Radon-Nikodym derivative of $\nu$ relative 
to $\mu$. The entropy $H(\nu|\mu)$, always a positive quantity,
is the optimal constant that makes the entropy inequality
\begin{equation} \label{entropy}
E_{\nu}[f] \leq \frac{1}{\alpha} \bigg\{ H(\nu|\mu)+ \log E_{\mu}[e^{\alpha f}] \bigg\}.
\end{equation}
true for any bounded, measurable function $f$ and $\alpha>0$.
A trivial consequence of the entropy inequality (\ref{entropy}) helps 
us to estimate $\nu(A)$, where $A$ is some event,
\begin{equation} \label{tigh1}
\nu (A) \leq \frac{\log(2)+H(\nu|\mu)}{\log(1+\frac{1}{\mu (A) })}
\end{equation} 

In our context we use the entropy inequality for the distribution of 
the process started in nonequilibrium and the distribution of the process started
in the equilibrium. 
We have,
\begin{equation} \nonumber
E^{\mathrm{neq}} \bigg[  \bigg| \int_0^T  g(x(t))dt \bigg| \bigg] =
E^{\mathrm{neq}} \bigg[ \frac{1}{N\alpha } \bigg| \int_0^T N\alpha g(x(t))dt \bigg| \bigg] \leq
\end{equation}
\begin{equation} \nonumber
\leq \frac{1}{N\alpha }H( \nu_N^{\mathrm{neq}}|\nu_N^{\mathrm{eq}}) + \frac{1}{N\alpha } \log
E^{\mathrm{eq}} \bigg[ \exp \bigg( \alpha \bigg| \int_0^T Ng(x(t)) dt \bigg| \bigg) \bigg] \leq
\end{equation} 
\begin{equation} \nonumber
\leq \frac{C}{\alpha} + \frac{1}{N\alpha } \log E^{\mathrm{eq}} \bigg[ 
\exp \bigg( N\alpha \int_0^T g(x(t)) dt \bigg) 
+ \exp \bigg( - N\alpha \int_0^T g(x(t)) dt \bigg) \bigg].
\end{equation}
As a consequence of the inequality 
$$\log(a+b) \leq \max (\log(2a), \log(2b)) \leq \log(2) + \max (\log(a), \log(b))$$  
for two positive numbers $a$ and $b$,  (\ref{l1}) follows as soon as we have
$$\limsup_{N \to \infty} \frac{1}{N}\log E^{\mathrm{eq}} \bigg[ \exp \bigg( N\alpha \int_0^T g(x(t)) dt \bigg) \bigg] \leq 0, \;  \; \alpha \neq 0.$$
A trivial consequence of  Feynman-Kac proves our Lemma,
\begin{equation} \label{feynman}
\frac{1}{N}\log E^{\mathrm{eq}} \bigg[ \exp \bigg( N\alpha \int_0^T g(x(t)) dt \bigg) \bigg] 
\leq T \mathrm{supspec}_{L^2(\nu^{\mathrm{eq}}_N )}(\alpha g +\frac{N^4}{N}L_N). 
\end{equation}
\qed 

{\bf The microscopic current $w$ can be replaced by 
local average of currents (the term (4.3)).} We shall show that the term (4.3) 
converges to zero in probability. As a consequence the current $w$ 
is replaced by a local average of $w$, that is closer to a deterministic
value. $w$ by itself is a single fluctuating random variable. 

We check that the hypothesis of Lemma \ref{lem1}  is valid 
for the function 
$g_{l,N}=\frac{1}{N} \sum_{i=1}^N \phi(\fri )N^2 ( w_i -  \Av{w}{i}{l} ) $. 
Recall that if the test function $\phi$ has continuous second-order derivative 
the quantity
$$\Avph -\phi(\fri ) = \frac{\phi(\frac{i-l}{N})-\phi(\fri)+ \dots +\phi(\frac{i+l}{N})-\phi(\fri)}{2l+1}$$
is of order ${\mathcal{O}}(\frac{l^2}{N^2})$.
Now, on integrating by parts twice, it follows for a fixed function $\rho$, with $E^{\mathrm{eq}}[\rho^2]=1$, that,
\begin{equation} \nonumber
\bigg| E^{\mathrm{eq}} \bigg[ \rho^2 \frac{1}{N}\sum_{i=1}^N \phi(\fri )N^2 ( w_i -  \Av{w}{i}{l} ) \bigg] \bigg| ^2 = 
\end{equation}
\begin{equation} \nonumber
=4 \bigg| E^{\mathrm{eq}} \bigg[ \frac{N^2}{N}\sum_{i=1}^N \sqrt{a_i}X_i(\rho) \sqrt{a_i}\Big(\Avph -\phi(\fri ) \Big) \rho  \bigg]   
\bigg|^2 \leq 
\end{equation}
\begin{equation} \nonumber
=4 E^{\mathrm{eq}} \bigg[ \sum_{i=1}^N \frac{N^4}{N^2} a_i \Big(\Avph -\phi(\fri ) \Big)^2 \rho ^2 \bigg] D_N(\rho) \leq   
C \frac{l^4}{N^4} \frac{N^4}{N} D_N(\rho) .
\end{equation} 
Therefore,
$$ \limsup_{l, N \to \infty} \sup_{\rho, E^{\mathrm{eq}}[\rho^2]=1}
   \frac{
\Big(E^{\mathrm{eq}} \Big[ \rho^2 \frac{N^2}{N}\sum_{i=1}^N \Big(\Avph -\phi(\fri ) \Big) w_i \Big]
\Big)^2  }{D_N(\rho)} \frac{N}{N^4} =0, $$
and hence (\ref{sup1}) is satisfied.  Moreover, 
\begin{equation} \nonumber
\lim_{l, N \rightarrow \infty} P^{\mathrm{neq}}_{N,T}\bigg\{ \bigg| \int_0^T \frac{1}{N} 
\sum_{i=1}^N \phi(\fri )N^2 ( w_i -  \Av{w}{i}{l}) dt \bigg| \geq \epsilon \bigg\} =0,  \; 
 \; \e >0.
\end{equation}

{\bf Inserting the fluctuations (the term (4.5)).} Let  $f(x_{-s},\dots,x_{s},\bar{x}_{0,l}) \in C^2(\mathbb{R}^{2s+1})$ be a local function that depends on the slope 
configuration in a box of size $s$ and on the mean 
slope in a box of a large size $l$. We use the notation
$l_1= l- \sqrt{l}$ for a slightly smaller $l$.
We want to show that the fluctuations approach zero, in the limit or that
\begin{equation} 
\lim_{l,N \to \infty} P^{\mathrm{neq}}_{N,T}\bigg\{ \bigg| \int_0^T \frac{N^2}{N} 
\sum_{i=1}^N \phi\bigg(\fri \bigg) L_N ( \Ave{f}{i}{l_1}) dt \bigg| \geq \epsilon \bigg\} =0,  \; \; \epsilon>0.
\end{equation}
We apply Ito formula, 
$$ \int_0^T \frac{N^2}{N} \sum_{i=1}^N \phi\bigg(\fri \bigg) L_N ( \Ave{f}{i}{l_1} ) dt  = $$
\begin{equation} \label{eq1}
=  \frac{1}{N\cdot N^2} 
\sum_{i=1}^N \phi\bigg(\fri \bigg) \bigg[ (\Ave{f}{i}{l_1})(x(T))-(\Ave{f}{i}{l_1} ) (x(0)) \bigg] + \frac{1}{N^2} M_N(t). 
\end{equation}

The summand in (\ref{eq1}) converges to zero, as the function $f$ is bounded. The second part of 
(\ref{eq1}) approaches zero because the $L^2$ norm of $\frac{M_N(T)}{N^2}$ is of order 
$\mathcal{O}(\frac{l^2}{N})$, as we can see below. 
Let $g=\frac{1}{N} \sum_{i=1}^N \phi(\fri ) (\Ave{f}{i}{l_1})(x)$, then  
$$ E^{\mathrm{neq}}\bigg[\frac{M_N(T)^2}{N^4}\bigg] = E^{\mathrm{neq}} \bigg[ \int_0^T L_N g -2gL_N g dt \bigg]  
= E^{\mathrm{neq}} \bigg[ \int_0^T \sum_{i=1}^N a_i (X_i g)^2 dt \bigg] = $$
$$ = \frac{1}{N^2} E^{\mathrm{neq}} \bigg[ \int_0^T \sum_{i=1}^N a_i \bigg(\sum_{k=1}^N \phi\bigg(\frac{k}{N}\bigg) X_i(\Ave{f}{k}{l_1}) \bigg) ^2 dt \bigg]. $$
The function $\Ave{f}{k}{l_1}$ is a cylinder function that depends just on the sites $n$ such that $k-l \leq n \leq k+l$. Therefore the vector field $X_i$ is zero when acting on most of the summands inside $\Ave{f}{k}{l_1}$. There are no more than  $2l$ sites $k$ such that $X_i(\Ave{f}{k}{l_1}) \neq 0$. We put all these arguments together to conclude that $ E^{\mathrm{neq}} [ M^2 (t)/N^4] \leq C l^2/N. $ \\

{\bf Replacing the current by the Laplacian of the slope field (the term (4.4)).}
In our model the instantaneous current, $w$, can not be written as  
 the discrete Laplacian $\tau h -2h -\tau^{-1}h$ of some local function $h$, thus
we use the method of Varadhan \cite{Var1} for computing the hydrodynamic scaling
limit of our model. 
The main idea is that the current decomposes as 
\begin{equation} \label{fluc}
 w = \hat{a}(\bar{x}_l) \Delta x + L_{\infty}f
\end{equation} 
for a suitable coefficient $\hat{a}(\bar{x}_l)$, where $\bar{x}_l$
is the average slope in a cube centered at the origin of microscopic side $l$.
The equation (\ref{fluc}) is known in the literature as the fluctuation-dissipation 
equations. We have explained before that terms of the form $L_{\infty}f$
have no effect on the macroscopic scale. 
A new feature is characteristic
to our model due to the complex interaction of the system: 
after filtering off the fluctuations from the current 
we are left with the Laplacian of some function and not 
a gradient, as happened for models previously considered 
in the literature. The precise meaning of (\ref{fluc}) is given below.

Our aim is to prove the existence of a sequence $\{f_r\}_{r \geq 0}$ of 
local functions and of the transport coefficient $\hat{a}$ such that 
\begin{eqnarray} \label{prob}
& &\limsup_{r,l,N \rightarrow \infty} P_{N,T}^{\mathrm{neq}} \Big\{ \Big| 
   \int_0^T \frac{N^2}{N}
   \sum_{i=1}^N \phi\bigg( \fri \bigg) \big[ \Av{w}{i}{l_1} - 
   \hat{a}(\bar{x}_{i,l}) \Av{(\Delta x)}{i}{l_1} - \nonumber \\
& & \quad \quad \quad \quad - \Ave{L_{\infty}f_r}{i}{l_1} \big] dt 
\Big| \geq \e \Big\} = 0.
\end{eqnarray}
The local function $f_r( x_{-s}, \dots , x_s, \bar{x}_{0,l})$ 
depends on two arguments, the mean slope $\bar{x}_{0,l}$ in a box of size $l$ and the 
slope field inside a box of size $s$, the size $s$ being much smaller than $l$. We can assume that 
the operator $L_{\infty}$ does not act on the first argument $\bar{x}_{0,l}$, 
since we can show that the action of the operator $L_{\infty}$ at the 
boundary sites is negligible. To be more precise, $\Ave{L_{\infty}f_r}{i}{l_1}$ stands for 
$$ \frac{L_{\infty}f_r(\bar{x}_{i,l}, \tau_{i-l_1}x) + \dots +L_{\infty}f_r(\bar{x}_{i,l}, 
\tau_{i+l_1}x)}{2l_1+1}. $$
Note that the function $\Ave{L_{\infty} f_r}{i}{l_1}$ depends just on the value of the field 
inside the box centered at $i$ and of size $l$.

As before we use Lemma \ref{lem1} to conclude that the event (\ref{prob}) 
has negligible probability in the limit. An additional difficulty shows up. 
If  $\lambda_{\e}$ is the largest eigenvalue of the 
perturbation, $L+\epsilon W $, of a negative operator $L$ with 
principal eigenvalue $0$, the eigenvalue $\lambda_{\e}$ 
has the formal series expansion,
$$\lambda_\epsilon= 0 +\epsilon E_{\nu}[W]+ \e^2 
<W, (-L)^{-1}W>_{\nu}+ {\mathcal O}(\epsilon^3),$$
Hence if the potential $W$ has mean zero, one expects 
$\lim_{\e \to 0} \lambda_{\e} \e^{-2}= <W, (-L)^{-1}W>_{\nu}$.
 Fortunately for suitable potential $W$ the central limit 
variance $<W, (-L)^{-1}W>_{\nu}$ converges to zero.

We shall need in our context, a particular result about the largest eigenvalue of a
perturbation operator, whose proof is found in Quastel \cite{Qua2}.

\begin{lemma} \label{clt2}
Let $W$ be a real potential that satisfies
$$<u, Wu>_{\nu} \leq l^{-1/2} D_l(u)^{1/2} ||u||_2,  \quad l < C \epsilon ^{-2/5}$$
for some $C$ small enough, or
$$||W||_{\infty} \leq C, \quad  \quad  K \leq (C \epsilon)^{-1/5}.$$
Provided that the generator $l^4L_l$ has spectral gap of order one the following estimate holds
\begin{equation}
\epsilon ^{-2} l^{-5} \mathrm{supspec}_{L^2{(\nu)}} (l^4 L_l - \epsilon l^5 W) \leq l <W, (-L_l)^{-1} W>_{\nu} + \mathcal{O}(1).
\end{equation}
 \end{lemma}


{\it Spectral gap.} Indeed the generator $l^4L_l$ of our model, 
defined by (\ref{genla}), has  a spectral gap of order $1$. 
The proof is standard by  the method  of Bakry-Emery (see Chang and Yau \cite{Cha} 
or Deuchel and Stroock \cite{Deu}).

The operator $L_l$ is an unbounded operator defined on the subspace $C_0^{\infty}(\RR^{2l+1})$
of the Hilbert space $L^2(\nu_{\alpha, l}^{gc})$ and is negative definite, with spectrum included 
in the negative semiaxis of the real line.  $0$ is an eigenvalue of the operator $L_{l}$ but 
the eigenspace corresponding to this eigenvalue is quite large, being infinite dimensional.

It is not hard to see that we can write the Hilbert space  $L^2 (\nu_{\alpha,l}^{gc})$  as
the direct sum $ \bigoplus_{y \in \RR ^2} L^2 ( \nu_{y,l}^{c} )$. Moreover 
because of the ergodicity of the dynamics (\ref{genla}) on the level sets of the function 
$y_{0,l}=(y^1_{0,l}, y^2_{0,l})$, we know that the eigenspace corresponding to zero 
of the restriction of the operator $L_{l}$  onto each Hilbert subspace $L^2 ( \nu_{y,l}^{c} )$,
is one dimensional. The next eigenvalue of the restriction  $L_{l}|_{L^2 (\nu_{y,l}^{c} )}$
is a negative number. The distance between the largest eigenvalue and the 
next largest eigenvalue of the operator $L_{l}|_{L^2 (\nu_{y,l}^{c} )}$ is called 
the spectral gap of the operator, because it is the gap in the spectrum of 
the operator.
\begin{lemma} \label{spec}(Spectral gap)
There is a universal constant $C$ that does not depend on 
the conserved quantities $y_{0,l}=(y_{0,l}^1,y_{0,l}^2) \in \RR^2$ such that
\begin{equation} \label{specgap}
\frac{C}{l^4}E_{\nu_{y,l}^c} [ \rho ^2 ] \leq  <(-L_{l})\rho, \rho>_{\nu_{y,N}^c}
\end{equation}
for any mean-zero function $\rho$, $E_{\nu_{y,l}^c}[\rho]=0$.
\end{lemma}

Since the generator $L_l$ of our model has a spectral gap, see Lemma \ref{spec},
 we can introduce the central limit theorem variance in our context.
\begin{definition} Suppose we have a cylinder function $f(x_{-s}, \dots x_s)$ such 
that  \\ $E_{\nu^{c}_{y,s}}[f] =0$ for all possible values of $y\in \RR^2$.
Recall that $\nu^{c}_{y,s}$ is the canonical measure in a box centered at 
$0$ and size $s$ defined in section \ref{model}. 
The central limit theorem variance of $f$ on the box $\Lambda_{i,l}$ 
is defined to be:
\begin{equation}
\VV_{i,l}(f, y) = 2(2l)<\Ave{f}{i}{l_1}, (-L_{i,l})^{-1} (\Ave{f}{i}{l_1})>_{\nu_{y,i,l}^c}.
\end{equation}  
If the box $\Lambda_{i,l}$ is centered at $0$ then we use the shorter notation 
$\VV_{l}(f, y)$ for the CLT-variance. At this point we stress that $\VV_{i,l}(f, y)$ 
is a local function depending on the field inside of $\Lambda_{i,l}$, more precisely
depending on the conserved quantities $y_{i,l}=(y_{i,l}^1, y_{i,l}^2)$, the mean slope 
and the linear mean of the slope field.
\end{definition}

At this point we stress that $\VV_{i,l}(f, y)$ 
is a local function depending on the field inside of $\Lambda_{i,l}$, more precisely
depending on the conserved quantities $y_{i,l}=(y_{i,l}^1, y_{i,l}^2)$, the mean slope 
and the linear mean of the slope field.

The strategy is to give a bound for the largest eigenvalue of 
a perturbation operator in terms of the CLT-variance $\VV_{i,l}(f, y)$.
Extra care must be taken because the CLT-variance $\VV_{i,l}(f, y)$ 
is uniformly-in-y small on bounded sets and not on unbounded sets.

The canonical measure $\nu^c_{i,l,y}$ for our model has been obtained by 
conditioning the grand canonical measure on the configurations 
with fixed mean slope and fixed linear mean slope in a box 
centered at $i$ and of size $l$. The second conditioning
makes the canonical measure not to have identical marginals.
Actually the expected values of the marginals depends linearly 
on the site. However the finite-dimensional marginals 
of the canonical measure converges toward the finite-dimensional
marginals of the grand canonical distribution as the size of the 
box approaches infinity, see Lemma \ref{equiven}. To benefit of this fact we will  
replace the CLT-variance $\VV_l(f,y)$ in a box of size $l$ with its expectation 
$E^{\mathrm{eq}}[\VV_l(f,y)|y_k]$ with respect to the canonical measure 
$\nu_{y,k}^c$ in a box of larger size $k$. We let $k$ go first to 
infinity. We formalize below. 

Let us define the function g as
\begin{equation}
g = \frac{N^2}{N} \sum_{i=1}^N \phi\bigg(\fri\bigg) \bigg[ \Av{w}{i}{l_1} - 
    \hat{a}(\overline{x}_{i,l})\Av{(\Delta x)}{i}{l_1} - 
    \Ave{L_{\infty}f_r}{i}{l_1} \bigg] 
\end{equation}
Thanks to Lemma \ref{lem1} the event (\ref{prob}) has negligible probability
if  $$\limsup_{r,l, N \to \infty} \mathrm{supspec}_{L^2(\nu_N^{\mathrm{eq}})}
(g + 2\beta \frac{N^4}{N}L_N)  \leq 0, \quad \quad \beta >0. $$
We write the operator $g + 2\beta \frac{N^4}{N}L_N$ as a sum of 
operators and we estimate the size of the principal eigenvalue 
of each operator in the sum. 
\begin{equation}
g + 2\beta \frac{N^4}{N}L_N = \Omega_1 + \Omega_2 +\Omega_3 +\Omega_4
\end{equation}
where,
\begin{eqnarray}
\quad \Omega_1 & = & g - \frac{1}{\beta N} \sum_{i=1}^N \phi\bigg(\fri\bigg)^2 
\VV_{i,l}(w-\hat{a}(\overline{x}_{i,l}) 
\Delta x -L_{\infty}f_r, y) 
+\beta \frac{N^4}{N}L_N  \nonumber \\
\quad \Omega_2 & = & \frac{1}{\beta N} \sum_{i=1}^N \phi\bigg(\fri\bigg)^2 
\bigg[\VV_{i,l}(w-\hat{a}(\overline{x}_{i,l}) 
\Delta x -L_{\infty}f_r, y)  
-  \nonumber \\
&  & -E^{\mathrm{eq}}[ \VV_{i,l}(w-\hat{a}(\overline{x}_{i,l}) \Delta x -
     L_{\infty}f_r, y) |y_{i,k}] \bigg] + \beta \frac{N^4}{N}L_N \nonumber \\
\quad \Omega_3 & = & \frac{1}{\beta N} \sum_{i=1}^N \phi\bigg(\fri\bigg)^2 
E^{\mathrm{eq}}[ \VV_{i,l}(w-\hat{a}(\overline{x}_{i,l}) \Delta x -
L_{\infty}f_r, y) |y_{i,k}] {\bf 1}_{|y_{i,k}|\geq \delta} \nonumber \\
\quad \Omega_4 & = &   \frac{1}{\beta N} \sum_{i=1}^N \phi\bigg(\fri\bigg)^2 
E^{\mathrm{eq}}[ \VV_{i,l}(w-\hat{a}(\overline{x}_{i,l}) \Delta x -
L_{\infty}f_r, y) |y_{i,k}] {\bf 1}_{|y_{i,k}|\leq \delta}. \nonumber
 \end{eqnarray}
The operators $\Omega_2$, $\Omega_3$ and $\Omega_4$ are 
understood as multiplication operators.

{\it The operator $\Omega_1$.}
Assume that $M$ is an upper bound for the test function $|\phi|$. Define 
$$g_{i,l} = \Av{w}{i}{l_1} - \hat{a}(\overline{x}_{i,l})\Av{(\Delta x)}{i}{l_1} 
- \Ave{L_{\infty}f_r}{i}{l_1}. $$
We have \\
$\mathrm{supspec}_{L^2(\nu_N^{\mathrm{eq}})} (\Omega_1) \leq$
\begin{eqnarray}
 & & \leq \sup_{|\lambda | \leq M} 
     \sup_{\rho, E^{\mathrm{eq}}[\rho ^2]=1} \bigg[ E^{\mathrm{eq}}\bigg[
     \lambda N^2 g_{0,l}\rho^2 -  
     \frac{\lambda ^2}{\beta}\VV_{l}(w-\hat{a}(\overline{x}_{l}) \Delta x -L_{\infty}f_r, y) 
     \rho ^2\bigg]-\nonumber \\ 
 & & \quad \quad \quad  -\frac{\beta N^4}{2l+1}D_l(\rho ) \bigg] \leq \nonumber \\ 
 & & \leq \sup_{|\lambda | \leq M} \sup_{y_l \in \RR ^2} \bigg[
\sup_{\rho, E^{\mathrm{eq}}[\rho ^2|y_l]=1} \bigg[ E^{\mathrm{eq}}[\lambda N^2 g_{0,l} \rho ^2] - 
\frac{\beta N^4}{2l+1}D_{\nu^c_{y,l}}(\rho) \bigg] - \nonumber \\
 & & \quad \quad \quad  
     -\frac{\lambda ^2}{\beta}\VV_{l}(w-\hat{a}(\overline{x}_{l}) \Delta x -L_{\infty}f_r, y)\bigg].
    \nonumber
\end{eqnarray}
Integrating by parts we can show that there is a constant $C$, not depending on $l$ and the values of the 
conserved quantities, such that for each density $\rho$:
$$E^{\mathrm{eq}}[g_{0,l} \rho^2|y_l] \leq
C\frac{\sqrt{D_{\nu^c_{y,l}}(\rho)}}{\sqrt{2l+1}}.$$
Hence the hypothesis of Lemma \ref{clt2} is satisfied and
$$ \limsup_{N \to \infty } \sup_{\rho, E^{\mathrm{eq}}[\rho ^2|y_l]=1} \bigg[ 
E^{\mathrm{eq}}[\lambda N^2 g_{0,l} \rho ^2] - 
\frac{\beta N^4}{2l+1}D_{\nu^c_{y,l}}(\rho ^2) \bigg] = $$
$$ =\frac{1}{(2l+1) l^4}\limsup_{N \to \infty } N^4 \sup_{\rho, E^{\mathrm{eq}}[\rho ^2|y_l]=1} 
\bigg[ E^{\mathrm{eq}}[\frac{\lambda(2l+1)l^4}{N^2} g_{0,l} \rho ^2] - 
\beta l^4 D_{\nu^c_{y,l}}( \rho ^2) \bigg]  \leq $$
$$ \leq \frac{\lambda ^2}{\beta} \VV_l(w-\hat{a}(\overline{x}_{l}) \Delta x -L_{\infty}f_r, y).$$
The convergence on the line above is uniform over the set of all possible values of the
conserved quantity $y_l= (y_l^1, y_l^2) \in \RR ^2$, therefore the principal eigenvalue
of the operator $\Omega_1$ becomes negative as $N \to \infty$.

{\it The operator $\Omega_2$.} 
Let us call 
$$v_{i,k}= V_{i,l}(w-\hat{a}(\overline{x}_{i,l}) \Delta x -L_{\infty}f_r, y) - 
E^{eq}[ V_{i,l}(w-\hat{a}(\overline{x}_{i,l}) \Delta x -L_{\infty}f_r, y) |y_{i,k}] .$$
We observe that $\Omega_2$ contains $v_{i,k}$ without being multiplied with
a factor $N^2$. Then,
$$\mathrm{supspec}_{L^2(\nu_N^{\mathrm{eq}})}(\Omega_2)=
\mathrm{supspec}_{L^2(\nu_N^{\mathrm{eq}})} \bigg[\frac{1}{\beta N} 
\sum_{i} \phi\bigg(\fri\bigg)^2 v_{i,k} + 
\beta \frac{N^4}{N}L_N \bigg] \leq $$
$$\leq \sup_{|\lambda | \leq M} \sup_{y_k \in \RR ^2} \bigg[
\sup_{\rho, E^{\mathrm{eq}}[\rho ^2|y_k]=1} \bigg[ E^{\mathrm{eq}}\bigg[
\frac{\lambda ^2}{\beta} v_{0,k} 
\rho ^2\bigg] - \frac{\beta N^4}{2k+1}D_{\nu^c_{y,k}}(\rho ) \bigg] = $$
$$= \frac{1}{(2k+1)k^4} \sup_{|\lambda | \leq M} \sup_{y_k \in \RR ^2} N^4 \bigg[\sup_{\rho, 
E^{\mathrm{eq}}[\rho ^2|y_k]=1} \bigg[ E^{\mathrm{eq}}\bigg[\frac{\lambda ^2(2k+1)k^4}{N^4\beta} v_{0,k} 
\rho ^2\bigg] - \beta k^4 D_{\nu^c_{y,k}}(\rho ) \bigg]. $$
Moreover, $v_{i,k}$ is a bounded function (see Lemma \ref{boundonvar}) and 
$E^{\mathrm{eq}}[v_{0,k}|y_k]=0$. We can apply Lemma \ref{clt2}.
$$ \limsup_{N \to \infty} N^4 \bigg[ \sup_{\rho, E^{\mathrm{eq}}[\rho ^2|y_k]=1} 
\bigg[ E^{\mathrm{eq}} \bigg[\frac{\lambda ^2(2k+1)k^4}{N^4\beta} v_{0,k} \rho ^2\bigg] 
- \beta k^4 D_{\nu^c_{y,k}}(\rho) \bigg] \leq 0 .$$  
The convergence on the line above is uniform over all possible values 
of the conserved quantities $y_l=(y_l^1, y_l^2)$; therefore $\limsup_{N \to \infty} \mathrm{supspec}_{L^2
(\nu_N^{\mathrm{eq}})}(\Omega_2) \leq 0.$ 

{\it The operator $\Omega_3$.} 
We refer to Lemma \ref{lem1}. Rather than proving that \\
$\limsup_{\delta, r, l, k, N \to \infty} \mathrm{supspec}_{L^2(\nu_N^{\mathrm{eq}})} (\Omega_3) =0$ we will show
\begin{equation} \label{wer*}
\limsup_{\delta, r, l, k, N \to \infty} \frac{1}{N} \log E^{\mathrm{eq}} 
\bigg[ \exp^{\int_0^T N\Omega_3} \bigg] ds =0,
\end{equation}
where
$$\Omega_3 = \frac{1}{\beta N}\sum_{i=1}^N \phi\bigg(\fri\bigg)^2 
E^{\mathrm{eq}}[V_l(w-\hat{a}(\bar{x}_l) \delta x -L_{\infty}f_r,y) | y_k] 
{\bf 1}_{|y_k| \geq \delta}.$$

The equation (\ref{wer*}) follows from the estimations
$$\frac{1}{N} \log E^{\mathrm{eq}}\bigg[ \exp \bigg(N \int_0^T \Omega_3 ds \bigg) \bigg] \leq
\frac{1}{N} \log E^{\mathrm{eq}}\bigg[ \frac{1}{T} \int_0^T \exp(N T \Omega_3) ds \bigg] =$$
$$= \frac{1}{N} \log E^{\mathrm{eq}}\bigg[ \exp \bigg(\sum_{i=1}^N \phi\bigg(\fri\bigg)^2 
E^{\mathrm{eq}}[\VV_l(w-\hat{a}(\bar{x}_l) \Delta x -L_{\infty}f_r,y) | y_k] 
{\bf 1}_{|y_k| \geq \delta}   \bigg) \bigg] \leq $$
$$ \leq \log E_{\nu^{gc}_{0}}[\exp (M^2 E^{\mathrm{eq}}
[\VV_l(w-\hat{a}(\bar{x}_l) \Delta x -L_{\infty}f_r,y) | y_k] {\bf 1}_{|y_k| \geq \delta}])].$$

The dominated convergence theorem can be applied (see the note at the end 
of section \ref{comp}) and gives us,
\begin{equation} \nonumber
 \lim_{\delta \mapsto 0} \log E_{\nu^{gc}_{0}}\Big[\exp (M^2 E^{\mathrm{eq}}
[\VV_l(w-\hat{a}(\bar{x}_l) \Delta x -L_{\infty}f_r,y) | y_k] {\bf 1}_{|y_k| \geq \delta}])\Big].
\end{equation}

{\it The operator $\Omega_4$.}
The main purpose of sections \ref{comp} is to show that there exists 
a sequence of functions $\{f_r \}_{r \geq 0}$ such that for any $\delta >0$.
$$ \limsup_{r, l, k \to \infty} \sup_{|y_k| \leq \delta} E^{\mathrm{eq}}[\VV_l(w-\hat{a}(\bar{x}_l) 
   \Delta x -L_{\infty}f_r,y) | y_k] =0, \quad \mathrm{and}$$
$$ \sup_{r} \sup_{y_k \in \mathbb{R}^2} E^{\mathrm{eq}}[\VV_l(w-\hat{a}(\bar{x}_l) 
   \Delta x -L_{\infty}f_r,y) | y_k] \leq  \infty . $$ 

{\bf The term (4.6)} It follows that the term (4.6) converges to zero in probability from
the following two lemmas. We refer the reader to Bertini, Olla and Landim \cite{Ber} 
or Guo, Papanicolau and Varadhan \cite{Guo} for the proof. The proof uses mainly the entropy 
inequality (\ref{entropy}) and consequence of Feynman-Kac formula (\ref{feynman}).

\begin{lemma} \label{local} (Local ergodicity) 
Let $f$ be a cylinder function. Define $\widetilde{f}: \mathbb{R} \to \RR$ to be the function $\widetilde{f}(\alpha ) = E_{\nu^{gc}_\alpha}[f]$. Then for any $\delta > 0$, $\phi: \mathbb{T} \to \RR$ a smooth function, 
$$ \limsup_{k,N \to \infty} P_{N,T}^{\mathrm{neq}} \bigg\{ \int_0^t \bigg| 
\frac{1}{N} \sum_{i=1}^N \phi\bigg(\fri\bigg)(\tau^i f(x(s)) -\widetilde{f} 
(\mathrm{Av}_{j=i-k}^{i+k} x_j (s))) \bigg| ds \geq \delta \bigg\}=0.$$
\end{lemma}

\begin{lemma}(Two-Block estimate) 
For any continuous function $g : \mathbb{R} \to \mathbb{R}$, let
$$ F_{k,a,N}=\bigg\{ \int_0^t \mathrm{Av}_{i=1}^N \mathrm{Av}_{j=i-aN}^{i+aN} 
(g(\mathrm{Av}_{l=i-k}^{i+k} x_l(s))-g(\mathrm{Av}_{l=j-k}^{j+k}x_l(s)))^2 ds 
\geq \delta \bigg\}.$$
Then for any $\delta>0$,
$$\limsup_{ k \to \infty, a \to 0, N \to \infty} P_{N,T}^{\mathrm{neq}}\{ F_{k,a,N}\} =0. $$
\end{lemma}

\section{Computation of the central limit theorem variances} \label{comp}
In this section we compute the value of the limit:
$$\limsup_{l \to \infty} \sup_{|\alpha|\leq \delta }E_{\nu_{\alpha}^{gc}}[\VV_{l}(f, y)]$$
for a particular class of cylinder functions $f$ to be described later.

We shall need a result that relates canonical and grand canonical measures,
known as the equivalence of ensemble. The equivalence of ensemble says that 
asymptotically  the marginal in a fixed box of the canonical measure is 
the marginal of the grand canonical measure. 
\begin{lemma} \label{equiven}  
(Equivalence of ensemble) Let $f(x_{-s}, \dots , x_s)$ be a bounded, local function. Then, 
for any $\e >0$, there exist $N_1 \in \NN$ and $\delta_1>0$ 
such that
$$ |E_{\nu^c_{y,k}}[ f] - E_{\nu^{gc}_\alpha}[f]| \leq \e$$
as long as $k> N_1$, $l>N_1$ and $|y_k^1-\alpha|<\delta_1$, $y_k^2 \in \RR$. 
\end{lemma} 

We shall need a new notation, namely $X_0^*$ for the adjoint of the vector field 
$X_0=\partial_{-1}-2\partial_0+\partial_1$ with respect to the inner-product of 
the Hilbert space $L^2(a \: d\nu_N^{\mathrm{eq}})$. 
The adjoint is given by the formula,
$$X_0^*(h) = -X_0(ah) -(V'(x_{-1})-2V'(x_0)+V'(x_1))ah.$$ 
Note that the current $w$ is equal to $ X_0^*(2)$, and the slope Laplacian $\Delta x$ 
is equal to $ X_0^*(\frac{1}{a})$.
\begin{lemma} \label{boundonvar}
Let $h(x_{-s}, \dots , x_s)$ be a bounded, cylinder function such that $f=X_0^*(h) \in 
L^2(\nu_{\alpha, s}^{gc})$ for all $\alpha \in \mathbb{R}$. Then there exists  
$C(h)< \infty $ that depends just on the function $h$ such that
\begin{equation}
\sup_{l, y_{0,l} \in \mathbb{R}^2} \frac{1}{l} < (-L_l)^{-1}(\sum_{|j|\leq l- \sqrt{l}} \tau^j f), 
\sum_{|i| \leq l-\sqrt{l}} \tau^i f>_{\nu_{y,l}^c}   \leq C(h).
\end{equation}
In addition if the function $h$ is bounded away from zero, $h \geq C > 0$, we have the lower bound,
$C_1(h) > 0$ that depends just on $h$,
\begin{equation}
C_1(h) \leq \sup_{l, y_{0,l} \in \mathbb{R}^2} \frac{1}{l} 
< (-L_l)^{-1}(\sum_{|j|\leq l- \sqrt{l}} \tau^j f), 
\sum_{|i| \leq l-\sqrt{l}} \tau^i f>_{\nu_{y,l}^c} . 
\end{equation}

\end{lemma}

\noindent
{\bf Proof.}
Call $\VV_l(f,y)=\frac{1}{l} < (-L_l)^{-1}(\sum_{|j|\leq l- \sqrt{l}} \tau^jf), 
\sum_{|i| \leq l-\sqrt{l}} \tau^if>_{\nu_{y,l}^c}$.
We use the variational formula
\begin{equation} \label{upvar}
\VV_l(f,y)= \sup_{u} \frac{<u, \sum_{|j|\leq l- \sqrt{l}} \tau^jf>_{\nu_{y,l}^c}^2}
{lD_{\nu_{y,l}^c}(u)}.
\end{equation}
We have
$$<u, \sum_{|j|\leq l- \sqrt{l}} \tau^jf>_{\nu_{y,l}^c}^2 = 
<u, \sum_{|j|\leq l- \sqrt{l}} X_j^*(\tau^j h)>_{\nu_{y,l}^c}^2=
E_{\nu_{y,l}^c} \bigg[ \sum_{|j|\leq l - \sqrt{l}} a_jX_j(u) \tau^jh \bigg]^2 $$
$$ \leq D_{\nu_{y,l}^c}(u) E_{\nu_{y,l}^c} \bigg[ \sum_{|j|\leq l- \sqrt{l}} a_j 
(\tau^jh)^2 \bigg]\leq 2l C(h) D_{\nu_{y,l}^c}(u),$$
therefore the upper bound for the variance is established.

For the lower bound, we may chose a particular function $u$ in (\ref{upvar}), namely 
$u = \sum_{i=-l}^l \frac{i^2}{2} x_i$. Note that $X_i(u)=1$ for any $-l+1 \leq i \leq l-1$.
It follows,
$$\frac{<u, \sum_{|j|\leq l- \sqrt{l}} \tau^jf>_{\nu_{y,l}^c}^2}{lD_{\nu_{y,l}^c}(u)} =
\frac{E_{\nu_{y,l}^c} \bigg[ \sum_{|j|\leq l - \sqrt{l}} a_j \tau^jh \bigg]^2 }
{l E_{\nu_{y,l}^c} [\sum_{|j| \leq l-1} a_j] } \geq C_1(h). $$
We have just used that $a$, and $h$ are functions bounded away from zero,
and $a$ is a bounded function. 
\newline
\qed

\noindent
{\it Note.} Lemma \ref{boundonvar} is true for any local function 
$f = \sum_{i=-s}^s X_j^*(h_j)$, where the functions $h_j$ are bounded, local functions.
\begin{definition}
For a bounded, local functions $f$ such that $E_{\nu_{y,l}^c}[f]=0$, for all possible values of 
$y \in \mathbb{R}^2$,  we define the semi-norm:
$$\lr{f}^2_{\alpha}=\limsup_{l \to \infty, y_k^1 \to \alpha} \frac{1}{l} 
E_{\nu_{y,k}^c}\bigg[< (-L_l)^{-1}(\sum_{|i| \leq l-\sqrt{l}} \tau^i f), 
\sum_{|j| \leq l-\sqrt{l}} \tau^j f>_{\nu_{y,l}^c}\bigg]. $$
\end{definition}
\noindent
We saw in Lemma \ref{boundonvar} that $\lr{f}_{\alpha}$ is a finite number as long as $f$ is equal 
to $X_0^*(h)$, where $h$ is a bounded, local function.
By polarization we can extend the semi-norm $\lr{ \quad }_{\alpha}$ to a semi-inner product 
$\lrp{\quad}{\quad }_{\alpha}$. 
For the remaining part of this section we compute the value of the semi-norm $\lr{f}_{\alpha}$ for  
certain functions $f$.
\begin{lemma}
Assume that $g$ is a bounded cylinder function with bounded first derivatives, 
$f =L_{\infty}g$, $w=X_0(a)-(V'(x_{-1})-2V'(x_0)+V'(x_1))a$ and 
$\Delta x =x_{-1}-2x_0+x_1$, then the following identities hold,
\begin{itemize}
\item[a)] $\lr{L_{\infty}g}_{\alpha}^2=
E_{\nu_{\alpha}^{gc}} \bigg[ a(x_{-1},x_0, x_1) \bigg( X_0 \big( \sum_{j \in \mathbb{Z}} \tau^j g \big) \bigg)^2 \bigg],$
\item[b)] $\lr{w}_{\alpha}^2=4E_{\nu_{\alpha}^{gc}}[a(x_{-1},x_0, x_1)],$
\item[c)] $\lrp{L_{\infty}g}{w}_{\alpha}= 2E_{\nu_{\alpha}^{gc}}[  a(x_{-1},x_0,x_1)
X_0\big(\sum_{j \in \mathbb{Z}} \tau^j g  \big)],$
\item[d)] $\lrp{L_{\infty}g}{\Delta x}_{\alpha}=0,$
\item[e)] $\lrp{w}{\Delta x}_{\alpha}=4 .$
\end{itemize}
\end{lemma}
\noindent
{\bf Proof.} 
One checks directly using equivalence of ensemble lemma \ref{equiven} and the asymptotic 
shift invariance of $\nu^c_{y,k}$ that the relations a)-e) hold.
In particular, for d) and e), it is important to notice that
$\Delta x = x_{-1}-2x_0+x_1 = X_0^*(\frac{1}{a})$. \\
\qed

\begin{definition}
Define the Hilbert space $\mathcal{H}_{\alpha}$ to be the closed linear span 
in $L^2(a\; d \nu_{\alpha}^{gc})$
of the function $1$ and functions $\xi_g=X_0(\sum_{j \in \mathbb{Z}} \tau^j g)$, 
where $g$ is a bounded local function with bounded first derivatives. 

\end{definition}

It is not hard to see that if $f$ is equal to $X_0^*(h)$ then, 
$$\lrp{f}{L_{\infty}g}_{\alpha}=E_{\nu_{\alpha}^{gc}}\bigg[a Proj_{\mathcal{H}_{\alpha}}(2h)
 X_0 \bigg( \sum_{j \in \mathbb{Z}} \tau^j g \bigg)\bigg], \quad
 \lrp{f}{w}_{\alpha}= E_{\nu_{\alpha}^{gc}}[a Proj_{\mathcal{H}_{\alpha}}(2h)2].$$
On both lines above $Proj_{\mathcal{H}_{\alpha}}$ stands for the projection 
opertor in the subspace $\mathcal{H}_{\alpha}$.
We are left to calculate $\lr{f}_{\alpha}$ for a function $f=X_0^*(h)$. 
As we will show in the following lemma,  
$\lr{f}_{\alpha} =E_{\nu_{\alpha}^gc}[a (Proj_{\mathcal{H}_{\alpha}}(2h))^2]$.

\begin{lemma}
Suppose $f$ is equal to $X_0^*(h)$, where $h$ is a bounded cylinder function. The value of the semi-norm of $f$, 
$$\lr{f}_{\alpha}= E_{\nu_{\alpha}^{gc}}[a (Proj_{\mathcal{H}_{\alpha}}(2h))^2].$$
\end{lemma}
\noindent
{\bf Proof.}
Using  Cauchy-Schwartz inequality it follows that,
$$\lr{f}_{\alpha}=\liminf_{l \to \infty, y_k^1 \to \alpha} \frac{1}{l} 
E_{\nu_{y,k}^c}\bigg[< (-L_l)^{-1}(\sum_{|i| \leq l-\sqrt{l}} \tau^i f), \sum_{|j| \leq l-\sqrt{l}} \tau^j f>_{\nu_{y,l}^c}\bigg]   \geq $$
$$\geq E_{\nu_{\alpha}^{gc}}[a( Proj_{\mathcal{H}_{\alpha}}(2h))^2].$$

Consider $g= \sum_{|j| \leq l- \sqrt{l}} \tau^j f$, where $f = X_0^*(h)$ then
$$\frac{1}{l} 
< (-L_l)^{-1}(\sum_{|i| \leq l- \sqrt{l}} \tau^i f), \sum_{|j| \leq l-\sqrt{l}} 
\tau^j f>_{\nu_{y,l}^c} =
\sup_{\rho,D_{\nu^c_{y,l}}(\rho)}\frac{<\sum_{|j| \leq l-\sqrt{l}} \tau^j f>^2_{\nu_{y,l}^c}}{l} =
$$
$$= \frac{ <\rho_l, \sum_{|j| \leq l-\sqrt{l}} \tau^j f>_{\nu_{y,l}^c} ^2}{2l^2} 
 = \frac{< \rho_l, \sum_{|j| \leq l- \sqrt{l}} X^*_j (\tau^j h)>_{\nu_{y,l}^c}^2}{2l^2}= $$
$$ = \frac{1}{ 2l^2} E_{\nu_{y,l}^c} \bigg[ \sum_{|j| \leq l-\sqrt{l}} a(x_{-j},x_j,x_j) X_j (\rho_l) \tau^j h \bigg]^2.$$
Above $\rho_l$ is some maximizer function with $D_{\nu^c_{y,l}}(\rho_l)= 2l$.

Call $u_l = \frac{1}{2l} \sum_{|j| \leq l-\sqrt{l}} \tau^{-j} (X_j \rho_l)$.
$$E_{\nu_{\alpha}^{gc}} [ a (u_l)^2] \leq 
  E_{\nu_{\alpha}^{gc}} \bigg[\frac{a}{2l} \sum_{|j| \leq l-\sqrt{l}}( \tau^{-j} 
   (X_j \rho_l))^2 \bigg] = \frac{1}{l} D_{\nu_{\alpha}^{gc}} (\rho_l) \leq 2.$$
\noindent
The above inequality shows that $\{u_l\}_l$ is a bounded sequence in $L^2(a  \, d \nu_{\alpha}^{gc} )$ 
and hence has a weakly convergent subsequence in $L^2(a \, d\nu_{\alpha}^{gc})$. 
Let $u$ be the weak limit of a convergent subsequence of $\{ u_l\}_l$. 
It follows that
$$ \limsup_{l \to \infty, y_k^1 \to \alpha} \frac{1}{l} 
E_{\nu_{y,k}^c}\bigg[< (-L_l)^{-1}(\sum_{|i| \leq l-\sqrt{l}} \tau^i f), \sum_{|j| \leq l-\sqrt{l}} \tau^j f>_{\nu_{y,l}^c}\bigg]   \leq $$
$$ \leq 2E_{\nu_{\alpha}^{gc}} [ a uh]^2 = 2 E_{\nu_{\alpha}^{gc}} [ a u Proj_{{\mathcal H}_{\alpha}}h]^2 \leq $$
$$\leq 2E_{\nu_{\alpha}^{gc}} [au^2] E_{\nu_{\alpha}^{gc}} [a(Proj_{{\mathcal H}_{\alpha}} h)^2] \leq  
E_{\nu_{\alpha}^{gc}} [a (Proj_{{\mathcal H}_{\alpha}} 2h)^2 ].$$

The key point that has allowed us to write the above inequalities is 
that the function $u$ has the property 
$X_a(\tau^b u) = X_b( \tau^a u)$ for all integers $a$ and $b$, and a 
function with this property belongs to $\mathcal{H}_{\alpha}$, 
(see Lemma \ref{clos1} or Savu \cite{Sav}, \cite{Sav1}). 
The next estimate shows that $X_a(\tau^b u) = X_b( \tau^a u)$ is valid
in a weak sense. Consider a smooth test function $\phi$, with bounded first 
derivatives. Assume $a>b$. We have
\begin{equation}
<X_a(\tau^b u_l), \phi>_{\nu^{gc}_{\alpha}}
 = < \frac{1}{2l} \sum_{|j| \leq l - \sqrt{l}} \tau^{b-j}(\rho_l), 
 X_b^* X_a^* \phi>_{\nu^{gc}_{\alpha}},\nonumber    
\end{equation}
and \\
$< X_a(\tau^b u_l) - X_b(\tau^a u_l) , \phi>_{\nu^{gc}_{\alpha}}=$ \\
\begin{eqnarray}
&&= < \frac{1}{2l} \sum_{j=b+1+l-\sqrt{l}}^{a+l - \sqrt{l}} \tau^{j}(\rho_l)
   - \frac{1}{2l} \sum_{j=b-l+\sqrt{l}}^{a-1-l + \sqrt{l}} \tau^{j}(\rho_l),
    X_b^* X_a^* \phi>_{\nu^{gc}_{\alpha}} = \nonumber \\
&&= < \frac{1}{2l} \sum_{j=b+1+l-\sqrt{l}}^{a+l - \sqrt{l}} \tau^{j}X_{b-j}(\rho_l)
   + \frac{1}{2l} \sum_{j=b-l+\sqrt{l}}^{a-1-l + \sqrt{l}} \tau^{j}X_{b-j}(\rho_l),
    X_a^* \phi>_{\nu^{gc}_{\alpha}}. \nonumber
\end{eqnarray}
By Cauchy-Schwartz inequality, we obtain \\
$< X_a(\tau^b u_l) - X_b(\tau^a u_l) , \phi>_{\nu^{gc}_{\alpha}}^2 \leq $
\begin{eqnarray} \label{covar}
& &\leq E_{\nu^{gc}_{\alpha}} \bigg[ \bigg( \frac{1}{2l} \sum_{j=b+1+l-\sqrt{l}}^{a+l - \sqrt{l}} \tau^{j}X_{b-j}(\rho_l)
   + \frac{1}{2l} \sum_{j=b-l+\sqrt{l}}^{a-1-l + \sqrt{l}} \tau^{j}X_{b-j}(\rho_l) \bigg)^2 \bigg]E_{\nu^{gc}_{\alpha}}[(X_a^* \phi)^2]\nonumber \\
&&\leq \frac{C(a-b)}{(2l)^2} \bigg( \sum_{j=b-a-l+\sqrt{l}}^{-1-l+\sqrt{l}}E_{\nu^{gc}_{\alpha}}[ (X_{j} \rho_l)^2]+
    \sum_{j=b-a+1+l-\sqrt{l}}^{l - \sqrt{l}} E_{\nu^{gc}_{\alpha}}[ (X_{j}\rho_l)^2] \bigg) \leq \nonumber \\
&&\quad \quad \quad \leq \frac{C(a-b)}{(2l)^2}D_{\nu^{gc}_{\alpha}}(\rho_l) \leq \frac{C(a-b)}{2l}. 
\end{eqnarray}
As $l$ converges to infinity, the sequence $\{ u_l \}_l$ approaches $u$ in the 
weak sense. Combining this fact with (\ref{covar}) we can establish $X_a(\tau^b u) = X_b( \tau^a u)$. 

Now we can conclude that if we have a local 
function $f$  such that $ f= X_0^*(h)$ then 
$$\lrp{f}{f}_{\alpha}  = \lim_{l \to \infty, y_k^1 \to \alpha} \frac{1}{l} 
E_{\nu_{y,k}^c}\bigg[< (-L_l)^{-1}(\sum_{|i| \leq l-\sqrt{l}} \tau^i f), 
  \sum_{|j| \leq l-\sqrt{l}} \tau^j f>_{\nu_{y,l}^c}\bigg]  =$$
$$=E_{\nu_{\alpha}^{gc}} [ a (Proj_{{\mathcal H}_{\alpha}}(2h))^2 ].$$
\qed

The results proved so far in this section allow us to conclude that if $g$ is a local function 
equal to $X_0^*(h)$ and $b(y_k^1)$ is a coefficient that depends on 
the mean-slope in a box of size $k$ then
$$\lrp{g + b w +L_{\infty} f}{g+ b w +L_{\infty} f}_{\alpha} = 
  \lim_{l \to \infty, y_k^1 \to \alpha}  
  E_{\nu_{y,k}^c}[ \VV_l(g+ b w +L_{\infty} f,y)]  = $$
$$=E_{a\nu_{\alpha}^{gc}} [(Proj_{{\mathcal H}_{\alpha}} (2h) +2b+ X_0(\sum_{i \in \mathbb{Z}} 
  \tau^i f))^2 ].$$
 
Since $Proj_{{\mathcal H}_{\alpha}}(2h)$ is a function in the closed linear span of 
the functions $1$ and $X_0(\sum_{i \in \mathbb{Z}} \tau^i f)$ where $f$ is a local
function it follows the existence of a sequence of functions $\{f_r\}_{r>0}$ and 
of a coefficient $b=\frac{1}{\hat{a}}$ such that 
$$ \lim_{r \to \infty} \lr{\Delta x +bw +L_{\infty} f_r}_{\alpha} =0. $$ 

\begin{lemma} \label{technical}
Let $a^* > 0$ be an upper bound for the function $a$. For any $\e > 0$ and $\delta > 0$ there exists a 
smooth function $f(x_{-l}, \dots ,x_l, \alpha)$ such that 
$$ \sup_{|\alpha| \leq \delta} \lr{w- L_{\infty}f- \hat{a}(\alpha)\Delta x}_{\alpha} 
< \e, \mathrm{and}$$
$$ \sup_{\alpha} \lr{w- L_{\infty}f- \hat{a}(\alpha)\Delta x}_{\alpha} <2a^* . $$
\end{lemma}

\noindent 
{\bf Proof.} We calculate 
$$\lr{w+ L_{\infty}f- \hat{a}(\alpha)\Delta x}_{\alpha} = 
4 \Big[ E_{\nu^{\mathrm{gc}}_{\alpha}}[a(1+X_0(\sum_{i \in \ZZ} \tau^if))^2]-\hat{a}(\alpha)\Big].$$
Let us introduce the notation $\mathcal{F}_{l,B}$ for the set of cylinder functions 
$f(x_{-l}, \dots , x_l)$ with $||f||_{\infty} \leq B$ and $||\partial_i f||_{\infty} \leq B$
for $i = -l, \dots , l$. Let $A(f, \alpha )= E_{\nu^{\mathrm{gc}}_{\alpha}}
[a(1+X_0(\sum_{i \in \ZZ} \tau^if))^2]-\hat{a}(\alpha)$ and $A_{l,B}(\alpha) = 
\inf_{f \in \mathcal{F}_{l,B}} A(f, \alpha)$.
The function $A_{l,B}(\alpha)$ is upper semicontinuous and nonincreasing in $l$ and $B$, and 
for each $\alpha \in \RR$, $\lim_{l,B \to \infty} A_{l,B}(\alpha)=0$, therefore,
 $$\lim_{l,B \to \infty} \sup_{|\alpha| \leq \delta +1} A_{l,B}(\alpha) =0.$$ 
Therefore we can find for each $\alpha$ a cylinder function $f(\alpha)$ such that 
$A(f(\alpha), \alpha) \leq \e$, if $|\alpha|\leq \delta +1$. We extend $f$ to be 
zero on $|\alpha|>\delta+1$. Then on $|\alpha|>\delta+1$ we have that 
$A(f(\alpha),\alpha ) = E_{\nu^{\mathrm{gc}}_{\alpha}} [a(x_{-1},x_0, x_1)] - \hat{a}(\alpha)$
and $A(f(\alpha),\alpha ) \leq a^* $. The avoid the problem that  $f$ is not smooth, we 
take the convolution of $f$ with a smoothing kernel $\phi$.  The required 
function is the convolution $f * \phi$. For a complete argument see Lemma 2.6 
in Quastel \cite{Qua6}.

\qed

\noindent
{\it  Note.} From the equivalence of ensemble lemma \ref{equiven}, 
we know that for any bounded cylinder function $f$ and any $\e >0$, 
there exit $N_1 \in \NN$ and $\delta_1>0$ such that
$$| E^{\mathrm{eq}}[\VV_l(w-\hat{a}(\bar{x}_{i,l}) \Delta x - L_{\infty}f, y)|y_k] -
 \lr{w - \hat{a}(\alpha ) \Delta x - L_{\infty}f}_{\alpha} | \leq \e,$$
as long as $k>N_1$, $l>N_1$, and $|y_k^1-\alpha|<\delta_1$.  
Then Lemma \ref{technical} helps us to conclude that for any $\e>0$ 
there exists a bounded cylinder function $f$ such that 
$$ \limsup_{l, k \to \infty} \sup_{|y_k| \leq \delta} E^{\mathrm{eq}}
[\VV_l(w-\hat{a}(\bar{x}_l) \Delta x -L_{\infty}f,y) | y_k] < 2\e, \quad \quad \mathrm{and}$$
$$ \sup_{y_k \in \mathbb{R}^2} E^{\mathrm{eq}}[\VV_l(w-\hat{a}(\bar{x}_l) \Delta x -L_{\infty}f_r,y) | y_k] \leq 3a^* . $$ 
$a^*$ on the line above is the upper bound for the function $a$ in the hypothesis of 
Lemma \ref{technical}.

{\bf Properties of the transport coefficient.}
The transport coefficient $\hat{a}$ has been defined as the unique real number such that
\begin{equation} \label{diff}
 \inf_f \lr{w - \hat{a}(\alpha) \Delta x - L_{\infty}f}^2_{\alpha} =0. 
\end{equation}
It is important to notice that the transport coefficient is also given by the formula:
$$ \hat{a}(\alpha) = \frac{\lrp{w}{\Delta x}_{\alpha}}{\lrp{\Delta x}{\Delta x}_{\alpha}} = 
   \frac{E_{\nu_{\alpha}^{gc}}[2\frac{2}{a}a]}{E_{\nu_{\alpha}^{gc}}[(Proj_{\mathcal{H}_{\alpha}}
   (\frac{2}{a}))^2a]}= \frac{1}{E_{\nu_{\alpha}^{gc}}[(Proj_{\mathcal{H}_{\alpha}} 
   (\frac{1}{a}))^2a]}=$$
$$ =\frac{1}{\sup_{\gamma, g} E_{\nu_{\alpha}^{gc}} [(Proj_{\gamma + X_0(\sum_{j \in \mathbb{Z}} 
\tau^jg)}\frac{1}{a})^2a]} = 
\inf_{\gamma,g} \frac{ E_{\nu_{\alpha}^{gc}}[(\gamma+ X_0(\sum_{j \in \mathbb{Z}} \tau^jg))^2a]}
{E_{\nu_{\alpha}^{gc}}[\frac{1}{a}(\gamma+X_0(\sum_{j \in \mathbb{Z}} \tau^jg))a]^2} =$$
\begin{equation} \label{diff1}
=\inf_{g} E_{\nu_{\alpha}^{gc}}\bigg[a(x_{-1},x_0,x_1)(1+ X_0(\sum_{j \in \mathbb{Z}}\tau^j g))^2\bigg].
\end{equation}

\section{Hydrodynamic limit of the model for surface electromigration}
\label{elect}

In this section we will show that the model for surface electromigration \ref{elect1}
has a hydrodynamic scaling limit as well. We will prove Theorem \ref{mainelect}.
The proof is inspired by Quastel \cite{Qua6}. 
From Cameron-Martin-Girsanov formula we can find the Radon-Nikodym
derivative 
\begin{eqnarray}
\frac{dP^{\mathrm{neq}}_{N,E,T}}{dP^{\mathrm{eq}}_{N,T}} & = &
\frac{d\nu_N^{\mathrm{neq}}}{d\nu_N^{\mathrm{eq}}}
\exp \bigg(
\frac{1}{2}\sum_{i=1}^N \int_0^T E\bigg(t, \frac{i}{N}\bigg) \sqrt{a(x_{i-1},x_i,x_{i+1})} dB_i - \nonumber \\
& &-\frac{1}{8} \sum_{i=1}^N \int_0^T E^2\bigg(t, \frac{i}{N}\bigg)a(x_{i-1},x_i,x_{i+1})dt \bigg). \nonumber
\end{eqnarray}
and the relative entropy
\begin{equation} \nonumber
H(P^{\mathrm{neq}}_{N,E,T}|P^{\mathrm{eq}}_{N,T}) = 
H(\nu_N^{\mathrm{neq}}|\nu_N^{\mathrm{eq}}) + 
E_{P^{\mathrm{neq}}_{N,E,T}} \bigg[\frac{1}{8} \sum_{i=1}^N \int_0^T E^2
\bigg(t, \frac{i}{N})a(x_{i-1},x_i,x_{i+1}\bigg)dt  \bigg].
\end{equation}
Hence there exists a constant $C$ such that 
$H(P^{\mathrm{neq}}_{N,E,T}|P^{\mathrm{eq}}_{N,T}) \leq CN$.
We can use the inequality (\ref{tigh1}) and the superexponential estimates 
(\ref{t1}), (\ref{t7}) to conclude that the sequence of probability measures
$\{Q_{N,E,T}^{\mathrm{neq}} \}_{N\geq 0}$ is tight in the topology of the 
space $\mathcal{X}$. Denote by $Q_{E,T}$ a 
weak limit of a subsequence of $\{Q_{N,E,T}^{\mathrm{neq}} \}_{N\geq 0}$.

We proceed further to identify the limit $Q_{E,T}$. To save space we ignore the time 
dependence of the test function $\phi$.\\
$\frac{1}{N}\sum_{i=1}^N \phi(\fri)x_i(T) - \frac{1}{N}\sum_{i=1}^N \phi(\fri)x_i(0) =$ 
\begin{eqnarray}
& &  = \int_0^T \frac{1}{2N}\sum_{i=1}^N \bigg( N^2 w_i+E\bigg(t, \fri\bigg)a_i \bigg) \phi''(\fri )dt + M_N(T) = \nonumber \\
& &  = \int_0^T \frac{N^2}{2N}\sum_{i=1}^N (w_i - \hat{a}(\bar{x}_{i,k})(\Delta x)_i - L_{\infty}\tau^i f_r) \phi''(\fri) dt + M_N(T)\nonumber \\
& & \quad  + \int_0^T \frac{N^2}{2N}\sum_{i=1}^N L_{N,E} \tau^i f_r \phi''(\fri) dt +   \\
& & \quad + \int_0^T \frac{1}{2N}\sum_{i=1}^N \bigg[ E\bigg(t, \fri \bigg)(a_i-\hat{a}(\bar{x}_{i,k}))+N^2(L_{\infty}\tau^if_r-L_{N,E} \tau^if_r) \bigg] \phi''(\fri) dt \\
& & \quad + \int_0^T \frac{1}{2N}\sum_{i=1}^N \hat{a}(\bar{x}_{i,k})\bigg[N^2(\Delta x)_i+E\bigg(t,\fri\bigg)\bigg] \phi''(\fri) dt \nonumber
\end{eqnarray}
Above $\{M_N(t)\}_{t \geq 0}$ is a martingale and as in section 4 we can 
prove that is negligible in the limit.
Because of the entropy inequality (\ref{entropy}) and the superexponential estimate
\begin{eqnarray} 
& &\limsup_{r,l,N \rightarrow \infty} \frac{1}{N} \log E^{\mathrm{eq}} \bigg[ \exp \Big( \Big| 
\int_0^T N^2\sum_{i=1}^N \big[ \Av{w}{i}{l_1} - \hat{a}(\bar{x}_{i,l}) \Av{(\Delta x)}{i}{l_1} - \nonumber \\
& & \quad \quad - \Ave{L_{\infty}f_r}{i}{l_1} \big] ds 
\Big| \Big) \bigg] \leq  0, \nonumber
\end{eqnarray}
the event on line (\ref{prob}) is negligible under $P^{\mathrm{neq}}_{N,E,T}$.
 For the new model that we have 
the negligible fluctuations are $L_{N,E}f$ and not $L_{\infty}f$ (i.e., 
the term (6.1) has no contribution towards the limit of the model).
The contribution coming from nontrivial fluctuations $L_{N,E}f-L_{\infty}f$ 
are gathered in the coefficient in front of the vector field $E$ in the 
nonlinear equation (\ref{noneqel}). Because of this reason, the coefficient turns 
out to be the transport coefficient $\hat{a}$, defined by the variational formula 
(\ref{coeff}) and not $\widetilde{a}(\alpha) = 
E_{\nu_{\alpha}^{gc}}[a]$    

We shall show that the sequence of smooth local functions $\{ f_r \}_{r\geq 0}$ 
introduced in the {\it Note} right after the Lemma \ref{technical} can be chosen 
to have the additional property that the term (6.2) is negligible in the
limit. Note that the term (6.2) is negligible if
$\limsup_{r,k,N \to \infty} P_{N,E,T}^{\mathrm{neq}}\{ |F(T)|> \e \} =0$
where,
\begin{equation}
F(T)= \int_0^T \frac{1}{2N}\sum_{i=1}^N E\bigg(t, \fri \bigg)\bigg[
a_i g_i - \hat{a}(\bar{x}_{i,k})\bigg] dt, \quad g_i = 1+X_i\bigg(
\sum_{j \in \ZZ} \tau^j f_r \bigg)
\end{equation}
We can write
\begin{eqnarray}
F(T) & = & \int_0^T \frac{1}{2N}\sum_{i=1}^N E\bigg(t, \fri \bigg)(a_i g_i 
    - E_{\nu^{gc}_{\bar{y}_{i,k}^1}}[a_i g_i]) + \\
& & +\int_0^T \frac{1}{2N}\sum_{i=1}^N E\bigg(t, \fri \bigg)(
    E_{\nu^{gc}_{\bar{y}_{i,k}^1}}[a_i g_i] - \hat{a}(\bar{x}_{i,k})) 
    {\bf 1}_{|\bar{y}_{i,k}^1| \leq \delta}+ \\
& & +\int_0^T \frac{1}{2N}\sum_{i=1}^N E\bigg(t, \fri \bigg)(
    E_{\nu^{gc}_{\bar{y}_{i,k}^1}}[a_i g_i] - \hat{a}(\bar{x}_{i,k})) 
    {\bf 1}_{|\bar{y}_{i,k}^1| \geq \delta} 
\end{eqnarray}
Recall that in the proof of Lemma \ref{technical} we introduced $A(f, \alpha)$ to be
the difference $E_{\nu^{\mathrm{gc}}_{\alpha}}
[a(1+X_0(\sum_{i \in \ZZ} \tau^if))^2]-\hat{a}(\alpha)$. 
We will show that we can modify the sequence $\{ f_r\}_{r\geq 0}$ to have the properties 
\begin{equation} \label{prop1}
\sup_{\alpha,r}A(f_r(\alpha),\alpha) \leq a^*, \quad 
\sup_{\alpha,r}|E_{\nu^{gc}_{\alpha}}[ag_0]-\hat{a}(\alpha) | \leq 5a^*,
\end{equation}
\begin{equation} \label{prop2} 
\lim_{r \to \infty} \sup_{|\alpha|\leq \delta} A(f_r(\alpha),\alpha) =0, \mathrm{and}
\end{equation}
\begin{equation} \label{prop3}
\lim_{r \to \infty} \sup_{|\alpha|\leq \delta} 
|E_{\nu^{gc}_{\alpha}}[ag_0]-\hat{a}(\alpha) | =0.
\end{equation}
We write,
\begin{equation} \nonumber
E_{\nu^{gc}_{\alpha}}[ag_0]-\hat{a}(\alpha) = 
(E_{\nu^{gc}_{\alpha}}[ag_0]- E_{\nu^{gc}_{\alpha}}[ag_0^2])+
(E_{\nu^{gc}_{\alpha}}[ag_0^2]-\hat{a}(\alpha))=
\end{equation}
\begin{equation} \nonumber
= E_{\nu^{gc}_{\alpha}}[ag_0(g_0-1)]+A(f_r(\alpha), \alpha).
\end{equation}
Since $|E_{\nu^{gc}_{\alpha}}[ag_0(g_0-1)]| \leq 4 a^*$ and 
$A(f_r(\alpha), \alpha) \leq a^*$ uniform in $\alpha$
and $r$, (\ref{prop1}) follows.
Consider  $f^*(\alpha)$ to be the minimizer of $A(f,\alpha)$ among 
local functions of $x_{-l}$ through $x_l$. $f^*$ has the property that for any $f_r$, 
$$E_{\nu^{gc}_{\alpha}}\bigg[ a\bigg(1+X_0(\sum_{j \in\ZZ}\tau^j f^*)\bigg)
X_0(\sum_{j \in\ZZ}\tau^j f_r)\bigg]=0.$$
Now, thanks to Cauchy-Schwartz inequality,
$$|E_{\nu^{gc}_{\alpha}}[ag_0(g_0-1)]| = 
\bigg|E_{\nu^{gc}_{\alpha}}\bigg[aX_0(\sum_{j \in \ZZ} \tau^j f_r) 
X_0(\sum_{j \in \ZZ}\tau^j(f^*-f_r))\bigg] \bigg|\leq$$
$$\leq E_{\nu^{gc}_{\alpha}}[a(X_0(\sum_{j \in \ZZ} \tau^j f_r))^2]^{1/2}
E_{\nu^{gc}_{\alpha}}[a(X_0(\sum_{j \in \ZZ}\tau^j(f^*-f_r)))^2]^{1/2}.$$
We modify $f_r$ such that
$$\sup_{|\alpha|\leq \delta}E_{\nu^{gc}_{\alpha}}\bigg[a\bigg(X_0(
\sum_{j \in \ZZ}\tau^j(f^*-f_r))\bigg)^2\bigg] \leq \frac{\e^2}{6a^*}.$$
To conclude (6.7) is negligible because of local ergodicity lemma \ref{local},
(6.8) is negligible because for any fixed $\delta$ we have (\ref{prop3}) and
(6.9) shrinks to zero once $\delta$ becomes very large. Also 
it is important to have property (\ref{prop1}). Our result follows. \\
\qed

\section*{Acknowledgement}
I would like to  thank Professor Jeremy Quastel for proposing me this
problem and for introducing me to interacting particle systems. 

\bibliographystyle{amsplain}

\begin{thebibliography}{11}

\bibitem{Bar} BARABASI, A. L. and STANLEY, H. E. {\it Fractal concepts in surface growth}, Cambridge University Press, Cambridge, 1995.
 
\bibitem{Ber} BERTINI, L., LANDIM, C. and OLLA, S. Derivation of 
Cahn-Hilliard Equations from Ginzburg-Landau Models. {\it J. of Statist. Physics} {\bf 88} (1997), no.1-2, 365-381.


\bibitem{Cha} CHANG, C. C. and YAU, H. T. Fluctuations of one-dimensional Ginzburg-Landau models in nonequilibrium.  {\it Comm. Math. Phys.}  {\bf 145} (1992), no. 2, 209--234.


\bibitem{Ede} EIDEL'MAN, S. D. {\it  Parabolic systems}. Amsterdam, North-Holland Pub. Co. 1969.

\bibitem{Deu} DEUSCHEL, J.-D. and STROOCK, D. {\it Large deviations}. Pure and 
Applied Mathematics, Academic Press, Boston, 1989.

\bibitem{Fri} FRITZ, J. On the hydrodynamic limit of a scalar Ginzburg-Landau lattice model: the resolvent approach.  {\it Hydrodynamic behavior and interacting particle systems.} IMA Vol. Math. Appl. {\bf 9} (1987), 75--97. 

\bibitem{Guo} GUO, M. Z., PAPANICOLAU, G. C. and  VARADHAN, S. R. S. Nonlinear diffusion limit for a system with nearest neighbor interactions.  {\it Comm. Math. Phys.}  {\bf 118} (1988), no. 1, 31--59.


\bibitem{Kip} KIPNIS, C. and LANDIM, C. {\it Scaling Limits of Interacting Particle Systems}. Springer-Verlag, Berlin, 1999.

\bibitem{Kru} KRUG, J., DOBBS, H. T. and MAJANIEMI, S. Adatom mobility for the solid-on-solid model. {\it Z. Phys. B} {\bf 97} (1995), 281--291.



\bibitem{Lan} LANDIM, C., OLLA, S. and VARADHAN, S. R. S. Symmetric simple exclusion
process: regularity of the self-diffusion coefficient. {\it Comm. Math. Phys.} 
{\bf 224} (2001), no. 1, 307--321.

\bibitem{Nis} NISHIKAWA, T. Hydrodynamic limit for the Ginzburg-Landau $\nabla\phi$ interface model with a conservation law. {\it J. Math. Sci. Univ. Tokyo}  {\bf 9} (2002), no. 3, 481--519. 

\bibitem{Qua1} QUASTEL, J. Diffusion of color in the simple exclusion process. {\it Comm. Pure Appl. Math.} {\bf XLV} (1992), 623-679.


\bibitem{Qua6} QUASTEL, J. Large deviations from a hydrodynamic 
scaling limit for a nongradient system. {\it Ann. Probab.} {\bf 23} (1995), no. 2, 724-742.
 


\bibitem{Sav} SAVU, A. {\it Hydrodynamic scaling limit of continuum solid-on-solid model.} Ph. D. thesis, 2004.
\bibitem{Sav1} SAVU, A. {\it Closed and exact functions in the context of Ginzburg-Landau models.} preprint.

\bibitem{Sch} SCHIMSCHAK, M. and KRUG, J. Surface electromigration as a moving boundary value problem. {\it Physical Review Letters} {\bf 78} (1997), 2, 278-281.

\bibitem{Var1} VARADHAN, S. R. S. Nonlinear diffusion limit for a system with nearest neighbour interactions II. {\it Asymptotic Problems in Probability Theory: Stochastic Models and Diffusion on Fractals (Sanda/Kyoto, 1990), 75-128, Pitman Research Notes in Mathematics}. {\bf 283} Longman Sci. Tech., Harlow, 1993.

\bibitem{Var2} VARADHAN, S. R. S. and  YAU, H. T. Diffusive limit of lattice gas with mixing conditions. {\it Asian J. Math.} {\bf 1} (1997), 4, 623-678.

\end{thebibliography}

\end{document}